\documentclass[11pt,a4paper,reqno]{amsart}
\usepackage{amssymb,euscript,graphicx,psfrag,amscd}

\newtheorem{theorem}{Theorem}

\DeclareMathOperator{\card}{card}
\DeclareMathOperator{\diam}{diam}
\DeclareMathOperator{\interior}{int}
\def\BB{\EuScript{B}}
\def\CC{\mathbb{C}}
\def\DD{\EuScript{D}}
\def\EE{\EuScript{E}}
\def\FF{\EuScript{F}}
\def\LL{\EuScript{L}}
\def\MM{\EuScript{M}}
\def\NN{\mathbb{N}}
\def\RR{\mathbb{R}}
\def\partition{\EuScript{R}}
\def\TT{\mathbb{T}}
\def\UU{\EuScript{U}}
\def\WW{\EuScript{W}}
\def\ZZ{\mathbb{Z}}

\begin{document}

\title[Hyperbolicity and Recurrence in Dynamical Systems]{Hyperbolicity
and Recurrence in Dynamical Systems: a survey of recent results}
\author{Luis Barreira}
\address{Departamento de Matem\'atica, Instituto
Superior T\'ecnico, 1049-001 Lisboa, Portugal}
\email{luis.barreira@math.ist.utl.pt}
\urladdr{http://www.math.ist.utl.pt/\textasciitilde barreira/}
\thanks{Partially supported by the Center for Mathematical Analysis,
Geometry, and Dynamical Systems, Lisbon, Portugal, through FCT's
Funding Program.}

\begin{abstract}
We discuss selected topics of current research interest in the
theory of dynamical systems, with emphasis on dimension theory,
multifractal analysis, and quantitative recurrence. The topics
include the quantitative versus the qualitative behavior of
Poincar\'e recurrence, the product structure of invariant measures
and return times, the dimension of invariant sets and invariant
measures, the complexity of the level sets of local quantities
from the point of view of Hausdorff dimension, and the conditional
variational principles as well as their applications to problems
in number theory.

We present the foundations of each area, and discuss recent
developments and applications. All the necessary notions from
ergodic theory, hyperbolic dynamics, dimension theory, and the
thermodynamic formalism are briefly recalled. We concentrate on
uniformly hyperbolic dynamics, although we also refer to
nonuniformly hyperbolic dynamics. Instead of always presenting the
most general results, we made a selection with the purpose of
illustrating the main ideas while we avoid the accessory
technicalities.
\end{abstract}


\maketitle
\tableofcontents

\vfill

\section*{Introduction}

Our main objective is to discuss selected topics of current
research interest in the theory of dynamical systems, with
emphasis on the study of \emph{recurrence}, \emph{hyperbolicity},
and \emph{dimension}. We recall that nontrivial recurrence and
hyperbolicity stand as principal mechanisms responsible for the
existence of stochastic behavior. We~want to proceed from
nontrivial recurrence and hyperbolicity, and discuss recent
results in three areas of research: \emph{dimension theory},
\emph{multifractal analysis}, and \emph{quantitative recurrence}.
We present a global view of the topics under discussion, although
the text substantially reflects a personal taste. In view of
readability, instead of always presenting the most general
results, we made a selection with the purpose of illustrating the
main ideas while we avoid the accessory technicalities.
Furthermore, all the necessary notions from hyperbolic dynamics,
ergodic theory, dimension theory, and the thermodynamic formalism
are briefly recalled. We~apologize if some reference was
overlooked, although if this happened it was of course totally
unintentional.

\subsection*{Recurrence}

The notion of nontrivial recurrence goes back to Poincar\'e in his
study of the three body problem. He proved in his celebrated
memoir~\cite{poincare} of 1890 that whenever a dynamical system
preserves volume almost all trajectories return arbitrarily close
to their initial position and they do this an infinite number of
times. This is Poincar\'e's recurrence theorem. The memoir is the
famous one that in its first version (printed in 1889, even having
circulated shortly, and of which some copies still exist today)
had the error that can be seen as the main cause for the study of
chaotic behavior in the theory of dynamical systems. Incidentally,
Poincar\'e's recurrence theorem was already present in the first
printed version of the memoir as then again in \cite{poincare}.
Already after the publication of~\cite{poincare}, the following
was observed by Poincar\'e about the complexity caused by the
existence of homoclinic points in the restricted three body
problem (as quoted for example in \cite[p.~162]{green}):
\begin{quote}
{\it One is struck by the complexity of this figure that I am not
even attempting to draw. Nothing can give us a better idea of the
complexity of the three-body problem and of all the problems of
dynamics in general~[...].}
\end{quote}
For a detailed and compelling historical account we recommend
\cite{green}.

\subsection*{Hyperbolicity}

The study of hyperbolicity goes back to the seminal work of
Ha\-da\-mard \cite{hada} in 1898 concerning the geodesic flow on
the unit tangent bundle of a surface with negative curvature, in
particular revealing the instability of the flow with respect to
the initial conditions. Hadamard observed (as quoted for example
in \cite[p.~209]{green}) that:
\begin{quote}
{\it [...] each stable trajectory can be transformed, by an
infinitely small variation in the initial conditions, into a
completely unstable trajectory extending to infinity, or, more
generally, into a trajectory of any of the types given in the
general discussion: for example, into a trajectory asymptotic to a
closed geodesic.}
\end{quote}
It should be noted that the geodesic flow preserves volume and as
such exhibits a nontrivial recurrence that was also exploited by
Hadamard. A~considerable activity took place during the 1920's and
1930's in particular with the important contributions of Hedlund
and Hopf who established several topological and ergodic
properties of geodesic flows, also in the case of manifolds with
not necessarily constant negative sectional curvature. We refer
the reader to the survey \cite{boris2} for details and further
references. Also in \cite{hada} Hadamard laid the foundations for
symbolic dynamics, subsequently developed by Morse and Hedlund and
raised to a subject in its own right (see in particular their
work~\cite{morse} of 1938; incidentally, this is where the
expression ``symbolic dynamics'' appeared for the first time).

\subsection*{Quantitative recurrence}

It should be noted that while the recurrence theorem of Poincar\'e
is a fundamental result in the theory of dynamical systems, on the
other hand it only provides information of qualitative nature. In
particular, it gives no information about the frequency with which
each trajectory visits a given set. This drawback was surpassed by
Birkhoff \cite{birk1, birk2} and von Neumann \cite{Neumann} who in
1931 established independently the first versions of the ergodic
theorem. Together with its variants and generalizations, the
ergodic theorem is a fundamental result in the theory of dynamical
systems and in particular in ergodic theory (one of the first
appearances of the expression ``ergodic theory'' occured in 1932
in joint work of Birkhoff and Koopman \cite{birk3}). Nevertheless,
the ergodic theorem considers only one aspect of the quantitative
behavior of recurrence. In particular, it gives no information
about the rate at which a given trajectory returns arbitrarily
close to itself. There has been a growing interest in the area
during the last decade, particularly with the work of Boshernitzan
\cite{bosher} and Ornstein and Weiss \cite{ow}.

\subsection*{Dimension theory}

In another direction, the dimension theory of dynamical systems
progressively developed, during the last two decades, into an
independent field. We emphasize that we are mostly concerned here
with the study of dynamical systems, and in particular of their
invariant sets and measures, from the point of view of dimension.
The first comprehensive reference that clearly took this point of
view is the book by Pesin~\cite{p}. The main objective of the
dimension theory of dynamical systems is to measure the complexity
from the dimensional point of view of the objects that remain
invariant under the dynamical system, such as the invariant sets
and measures. It turns out that the thermodynamic formalism
developed by Ruelle in his seminal work~\cite{R1} has a privileged
relation with the dimension theory of dynamical systems.

\subsection*{Multifractal analysis}

The multifractal analysis of dynamical systems is a subfield of
the dimension theory of dynamical systems. Briefly, multifractal
analysis studies the complexity of the level sets of invariant
local quantities obtained from a dynamical system. For example, we
can consider Birkhoff averages, Lyapunov exponents, pointwise
dimensions, and local entropies. These functions are usually only
measurable and thus their level sets are rarely manifolds. Hence,
in order to measure their complexity it is appropriate to use
quantities such as the topological entropy or the Hausdorff
dimension. Multifractal analysis has also a privileged relation
with the experimental study of dynamical systems. More precisely,
the so-called multifractal spectra obtained from the study of the
complexity of the level sets can be determined experimentally with
considerable precision. As such we may expect to be able to
recover some information about the dynamical system from the
information contained in the multifractal spectra.

\subsection*{Contents of the survey}

We first introduce in Section~\ref{sec1} the fundamental concepts
of ergodic theory and hyperbolic dynamics, and in particular
nontrivial recurrence and hyperbolicity. Section~\ref{sec2} is
dedicated to the discussion of the concept of nonuniform
hyperbolicity, that in particular includes the case of a dynamical
system preserving a measure with all Lyapunov exponents nonzero.
This theory (also called Pesin theory) is clearly recognized today
as a fundamental step in the study of stochastic behavior. We
consider hyperbolic measures in Section~\ref{sec3}, i.e., measures
with all Lyapunov exponents nonzero and describe their product
structure that imitates the product structure observed on
hyperbolic sets. We also discuss in Section~\ref{sec3} the
relation with the dimension theory of invariant measures.

Section~\ref{sec4} is dedicated to the study of the dimension
theory of invariant sets. This study presents complications of
different nature from those in the dimension theory of invariant
measures. In particular, the study of the dimension of invariant
sets is often affected by number-theoretical properties. Our
emphasis here is on the study of the so-called geometric
constructions which can be seen as models of invariant sets of
dynamical systems. There is again a privileged relation with the
thermodynamic formalism and we start to describe this relation in
Section~\ref{sec4}. Section~\ref{sec5} is dedicated to the study
of the dimension of invariant sets of hyperbolic dynamics, both
invertible and noninvertible. In particular, we present the
dimension formulas for repellers and hyperbolic sets in the case
of conformal dynamics. Symbolic dynamics plays a fundamental role
in some studies of dimension and is also considered in
Section~\ref{sec5}. In particular, this allows us to model
invariant sets by geometric constructions. Multifractal analysis
is the main theme of Section~\ref{sec6}. We describe the interplay
between local and global properties in the case of hyperbolic
dynamics. We also present several examples of multifractal spectra
and describe their properties.

In Section~\ref{sec7} we discuss the properties of the set of
points for which the Birkhoff averages do not converge. In view of
the ergodic theorem this set has zero measure with respect to any
invariant measure, and thus it is very small from the point of
view of measure theory. On the other hand, it is rather large from
the point of view of dimension theory and entropy theory. We also
discuss how one can make rigorous a certain multifractal
classification of dynamical systems. We discuss conditional
variational principles in Section~\ref{sec8}. These have a
privileged relation with multifractal analysis: roughly speaking,
the multifractal analysis of a given multifractal spectrum is
equivalent to the existence of a corresponding conditional
variational principle. We also discuss in Section~\ref{sec8} how
these variational principles and in particular their
multi-dimensional versions have applications to certain problems
in number theory. In Section~\ref{sec9} we address the problem of
quantitative recurrence allude to above. We also describe the
product structure of return times thus providing an additional
view of the product structure described before for invariant sets
and invariant measure. We also briefly describe some applications
to number theory.

\medskip
\noindent\emph{Acknowledgment.}  I am particularly indebted to my
friends and collaborators Yakov Pesin, Beno{\^\i}t Saussol, J\"org
Schmeling, and Christian Wolf with whom I enjoyed countless
mathematical as well as nonmathematical conversations.

\section{Hyperbolicity and nontrivial recurrence}\label{sec1}

\subsection{Dynamical systems and hyperbolicity}\label{sec1.1}

One of the paradigms of the theory of dynamical systems is that
the local instability of trajectories influences the global
behavior of the system and opens the way to the existence of
stochastic behavior. Mathematically, the instability of
trajectories corresponds to some degree of hyperbolicity.

Let $f\colon M\to M$ be a diffeomorphism and $\Lambda\subset M$ an
$f$-invariant set, i.e., a set such that $f^{-1}\Lambda=\Lambda$.  We
say that $\Lambda$ is a \emph{hyperbolic set} for $f$ if for every
point $x\in\Lambda$ there exists a decomposition of the tangent space
\begin{equation}\label{**1}
T_xM=E^s(x)\oplus E^u(x)
\end{equation}
varying continuously with $x$ that satisfies
\[
d_xfE^s(x)=E^s(fx) \quad\text{and}\quad d_xfE^u(x)=E^u(fx),
\]
and there exist constants $\lambda\in(0,1)$ and $c>0$ such that
\[
\lVert d_xf^n|E^s(x)\rVert\le c\lambda^n \quad\text{and}\quad
\lVert d_xf^{-n}|E^u(x)\rVert\le c\lambda^n
\]
for each $x\in\Lambda$ and $n\in\NN$.

Given $\varepsilon>0$, for each $x\in M$ we consider the sets
\[
V^s_\varepsilon(x)=\{y\in
B(x,\varepsilon):\text{$d(f^ny,f^nx)<\varepsilon$ for every
$n>0$}\}
\]
and
\[
V^u_\varepsilon(x)=\{y\in
B(x,\varepsilon):\text{$d(f^ny,f^nx)<\varepsilon$ for every
$n<0$}\},
\]
where $d$ is the distance on $M$ and $B(x,\varepsilon)\subset M$
is the open ball centered at~$x$ of radius~$\varepsilon$.

A hyperbolic set possesses a very rich structure.

\begin{theorem}[Hadamard--Perron Theorem]\label{thm:HP}
If $\Lambda$ is a compact hyperbolic set for a $C^1$
diffeomorphism then there exists $\varepsilon>0$ such that for
each $x\in\Lambda$ the sets $V^s_\varepsilon(x)$ and
$V^u_\varepsilon(x)$ are manifolds containing~$x$ and satisfying
\begin{equation}\label{**2}
T_xV^s_\varepsilon(x)=E^s(x) \quad\text{and}\quad
T_xV^u_\varepsilon(x)=E^u(x).
\end{equation}
\end{theorem}

We refer the reader to the book by Anosov \cite[\S4]{Anosov} for
references and for a detailed account of the origins of the
Hadamard--Perron Theorem.

\begin{figure}[htbp]
\begin{psfrags}
\psfrag{Vsx}{$V^s_\varepsilon(x)$}
\psfrag{Vsy}{$V^u_\varepsilon(x)$} \psfrag{Esx}{$E^s(x)$}
\psfrag{Eux}{$E^u(x)$} \psfrag{X}{$x$}
\begin{center}
\includegraphics{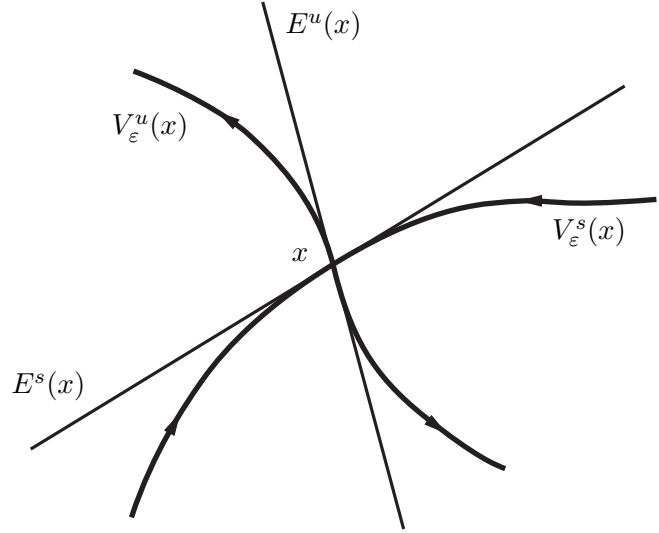}
\end{center}
\end{psfrags}
\caption{Local stable manifold and local unstable
manifold}\label{fig1}
\end{figure}

\begin{figure}[htbp]
\begin{psfrags}
\psfrag{Vsx}{$V^s_\varepsilon(x)$} \psfrag{X}{$x$} \psfrag{Y}{$y$}
\psfrag{Vsy}{$V^u_\varepsilon(y)$} \psfrag{Xy}{$[x,y]$}
\begin{center}
\includegraphics{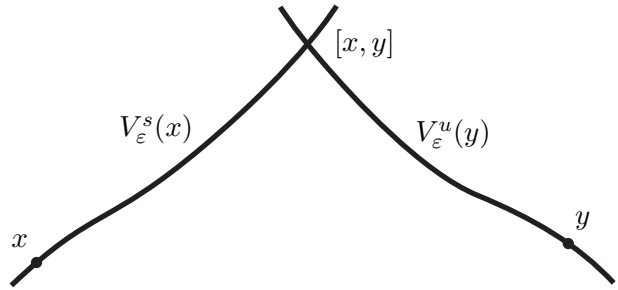}
\end{center}
\end{psfrags}
\caption{Product structure}\label{fig2}
\end{figure}

The manifolds $V^s_\varepsilon(x)$ and $V^u_\varepsilon(x)$ are
called respectively \emph{local stable manifold} and \emph{local
unstable manifold} at~$x$ (of size~$\varepsilon$). It follows from
\eqref{**1} and \eqref{**2} that these manifolds are transverse
(see Figure~\ref{fig1}). Furthermore, under the assumptions of
Theorem~\ref{thm:HP} one can show that the sizes of
$V^s_\varepsilon(x)$ and $V^u_\varepsilon(x)$ are uniformly
bounded away from zero, i.e., there exists
$\gamma=\gamma(\varepsilon)>0$ such that
\[
V^s_\varepsilon(x)\supset B^s_\gamma(x) \quad\text{and}\quad
V^u_\varepsilon(x)\supset B^u_\gamma(x)
\]
for every point $x\in\Lambda$, where $B^s_\gamma(x)$ and
$B^u_\gamma(x)$ are the open balls centered at~$x$ of radius
$\gamma$ with respect to the distances induced by $d$ respectively
on $V^s_\varepsilon(x)$ and $V^u_\varepsilon(x)$. The continuous
dependence of the spaces $E^s(x)$ and $E^u(x)$ in $x\in\Lambda$
guarantees that there exists $\delta=\delta(\varepsilon)>0$ such
that if $d(x,y)<\delta$ for two given points~$x$, $y\in\Lambda$
then the intersection $V^s_\varepsilon(x)\cap V^u_\varepsilon(y)$
is composed of exactly one point. We call \emph{product structure}
to the function
\[
[\cdot,\cdot]\colon\{(x,y)\in\Lambda\times\Lambda:d(x,y)<\delta\}\to M
\]
defined by $[x,y]=V^s_\varepsilon(x)\cap V^u_\varepsilon(y)$ (see
Figure~\ref{fig2}).

When the whole manifold $M$ is a hyperbolic set for $f$ we say
that $f$ is an \emph{Anosov diffeomorphism}. This class of
diffeomorphisms was introduced and studied by Anosov
in~\cite{Anosov}. The notion of hyperbolic set was introduced by
Smale in his seminal paper~\cite{smale}. Anosov diffeomorphisms
and more generally the diffeomorphisms with hyperbolic sets
constitute in a certain sense the class of transformations with
the ``strongest possible'' hyperbolicity. Moreover, hyperbolicity
is one of the main mechanisms responsible for the stochastic
behavior in natural phenomena, even though not always with the
presence of (uniformly) hyperbolic sets as defined above. These
considerations justify the search for a ``weaker'' concept of
hyperbolicity, present in a much more general class of dynamical
systems that we will call \emph{nonuniformly hyperbolic dynamical
systems} (see Section~\ref{sec2.1}).  The study of these systems
is much more delicate than the study of diffeomorphisms with
hyperbolic sets and namely of Anosov diffeomorphisms. However, it
is still possible to establish the presence of a very rich
structure and in particular the existence of families of stable
and unstable manifolds (see Section~\ref{sec2.1}).

\subsection{Ergodic theory and nontrivial recurrence}\label{sec1.2}

We now introduce the concept of invariant measure, which
constitutes another fundamental departure point for the study of
stochastic behavior.  Namely, the existence of a finite invariant
measure causes the existence of a \emph{nontrivial recurrence},
that is proper of stochastic behavior.

If $T\colon X\to X$ is a measurable transformation, we say that a
measure $\mu$ on~$X$ is \emph{$T$-invariant} if
\[
\mu(T^{-1}A)=\mu(A)
\]
for every measurable set $A\subset X$.  The study of
transformations with invariant measures is the main theme of
ergodic theory.

In order to describe rigorously the concept of nontrivial
recurrence, we recall one of the basic but fundamental results of
ergodic theory---the Poincar\'e recurrence theorem. This result
states that any dynamical system preserving a finite measure
exhibits a nontrivial recurrence in any set $A$ with positive
measure, in the sense that the orbit of almost every point in $A$
returns infinitely often to $A$.

\begin{theorem}[Poincar\'e recurrence theorem]\label{thm:pp}
Let $T\colon X\to X$ be a measurable transformation and $\mu$ a
$T$-invariant finite measure on~$X$.  If $A\subset X$ is a
measurable set with positive measure then
\[
\card\{n\in\NN:T^nx\in A\}=\infty
\]
for $\mu$-almost every point $x\in A$.
\end{theorem}

A slightly modified version of Theorem~\ref{thm:pp} was first
established by Poincar\'e in his seminal memoir on the three body
problem \cite{poincare} (see \cite{green} for a detailed
historical account).

The simultaneous existence of \emph{hyperbolicity} and
\emph{nontrivial recurrence} ensures the existence of a very rich
orbit structure. Roughly speaking, the nontrivial recurrence
allows us to conclude that there exist orbits that return as close
to themselves as desired. On the other hand, the existence of
stable and unstable manifolds at these points and their
transversality guarantees the existence of transverse homoclinic
points, thus causing an enormous complexity through the occurrence
of Smale horseshoes (see also Section~\ref{sec2}). As such,
hyperbolicity and nontrivial recurrence are two of the main
mechanisms responsible for the existence of stochastic behavior in
natural phenomena.

\section{Nonuniform hyperbolicity}\label{sec2}

\subsection{Nonuniformly hyperbolic trajectories}\label{sec2.1}

The concept of \emph{nonuniform hyperbolicity} originated in the
fundamental work of Pesin \cite{Pesin1, Pesin2, Pesin3}, in
particular with the study of smooth ergodic theory that today is
clearly recognized as a fundamental step in the study of
stochastic behavior (see \cite{lyap, km, M, viana} and the
references therein).

Let  $f\colon M\to M$ be a diffeomorphism.  The trajectory
$\{f^nx:n\in\ZZ\}$ of a point $x\in M$ is called
\emph{nonuniformly hyperbolic} if there exist decompositions
\[
T_{f^nx}M=E^s_{f^nx}\oplus E^u_{f^nx}
\]
for each $n\in\ZZ$, a constant $\lambda\in(0,1)$, and for each
sufficiently small $\varepsilon>0$ a positive function
$C_\varepsilon$ defined on the trajectory of~$x$ such that if
$k\in\ZZ$ then:
\begin{enumerate}
\item\label{ap2:xx} $C_\varepsilon(f^kx)\le e^{\varepsilon\lvert
k\rvert}C_\varepsilon(x)$; \item\label{ap2:2} $d_xf^kE^s_x =
E^s_{f^kx}$ and $d_xf^kE^u_x = E^u_{f^kx}$; \item\label{ap2:3} if
$v \in E^s_{f^k x}$ and $m>0$ then
\[
\lVert d_{f^kx}f^mv\rVert\le C_\varepsilon(f^kx)\lambda^m
e^{\varepsilon m}\lVert v\rVert;
\]
\item\label{ap2:4} if $v\in E^u_{f^k x}$ and $m<0$ then
\[
\lVert d_{f^kx}f^mv\rVert\le C_\varepsilon(f^kx)\lambda^{\lvert
m\rvert} e^{\varepsilon \lvert m\rvert}\lVert v\rVert;
\]
\item\label{ap2:yy} $\angle(E^u_{f^kx},E^s_{f^kx})\ge
C_\varepsilon(f^kx)^{-1}$.
\end{enumerate}

The expression ``nonuniform'' refers to the estimates in
conditions \ref{ap2:3} and~\ref{ap2:4}, that can differ from the
``uniform'' estimate $\lambda^m$ by multiplicative terms that may
grow along the orbit, although the exponential rate $\varepsilon$
in \ref{ap2:xx} is small when compared to the constant
$-\log\lambda$.  It is immediate that any trajectory in a
hyperbolic set is nonuniformly hyperbolic.

Among the most important properties due to nonuniform
hyperbolicity is the existence of stable and unstable manifolds
(with an appropriate version of Theorem~\ref{thm:HP}), and their
``absolute continuity'' established by Pesin in \cite{Pesin1}. The
theory also describes the ergodic properties of dynamical systems
with an invariant measure absolutely continuous with respect to
the volume \cite{Pesin2}. Also of importance is the Pesin entropy
formula for the Kolmogorov--Sinai entropy in terms of the Lyapunov
exponents \cite{Pesin2} (see also \cite{LY0}). Combining the
\emph{nonuniform hyperbolicity} with the \emph{nontrivial
recurrence} guaranteed by the existence of a finite invariant
measure (see Section~\ref{sec1.2}), the fundamental work of Katok
in \cite{Katok2} revealed a very rich and complicated orbit
structure (see also \cite{km}).

We now state the result concerning the existence of stable and
unstable manifolds, established by Pesin in \cite{Pesin1}.

\begin{theorem}[Existence of invariant manifolds]\label{thm:var}
If $\{f^nx:n\in\ZZ\}$ is a nonuniformly hyperbolic trajectory of a
$C^{1+\alpha}$ diffeomorphism, for some $\alpha>0$, then for each
sufficiently small $\varepsilon>0$ there exist manifolds $V^s(x)$
and $V^u(x)$ containing~$x$, and a function $D_\varepsilon$
defined on the trajectory of~$x$ such that:
\begin{enumerate}
\item $T_xV^s(x)=E^s_x$ and $T_xV^u(x)=E^u_x$; \item
$D_\varepsilon(f^kx)\le e^{2\varepsilon\lvert
k\rvert}D_\varepsilon(x)$ for each $k\in\ZZ$; \item if $y\in
V^s(x)$, $m>0$ and $k\in\ZZ$ then
\begin{equation}\label{9.1}
d(f^{m+k}x,f^{m+k}y)\le D_\varepsilon(f^kx)\lambda^m
e^{\varepsilon m}d(f^kx,f^ky);
\end{equation}
\item if $y\in V^u(x)$, $m<0$ and $k\in\ZZ$ then
\begin{equation}\label{9.1*}
d(f^{m+k}x,f^{m+k}y)\le D_\varepsilon(f^kx)\lambda^{\lvert
m\rvert} e^{\varepsilon \lvert m\rvert}d(f^kx,f^ky).
\end{equation}
\end{enumerate}
\end{theorem}

The manifolds $V^s(x)$ and $V^u(x)$ are called respectively
\emph{local stable manifold} and \emph{local unstable manifold}
at~$x$.  Contrarily to what happens with hyperbolic sets, in the
case of nonuniformly hyperbolic trajectories the ``size'' of these
manifolds may not be bounded from below along the orbit (although
they may decrease at most with an exponentially small speed when
compared to the speeds in \eqref{9.1} and~\eqref{9.1*}). This
makes their study more complicated.

In \cite{Pu}, Pugh constructed a $C^1$ diffeomorphism in a
manifold of dimension~$4$, that is not of class $C^{1+\alpha}$ for
any $\alpha>0$ and for which there exists no manifold tangent to
$E^s_x$ such that the inequality \eqref{9.1} is valid in some open
neighborhood of~$x$. This example shows that the hypothesis
$\alpha>0$ is crucial in Theorem~\ref{thm:var}.

The proof of Theorem~\ref{thm:var} in \cite{Pesin1} is an
elaboration of the classical work of Perron.  This approach was
extended by Katok and Strelcyn in \cite{KS} for maps with
singularities. In \cite{RR1}, Ruelle obtained a proof of
Theorem~\ref{thm:var} based on the study of perturbations of
products of matrices in the Multiplicative ergodic theorem (see
Theorem~\ref{Oset} in Section~\ref{sec2.2}).  Another proof of
Theorem~\ref{thm:var} is due to Pugh and Shub~\cite{PS} with an
elaboration of the classical work of Hadamard.  See \cite{lyap,
FHY, km} for detailed expositions.

There exist counterparts of Theorem~\ref{thm:var} for dynamical
systems in infinite dimensional spaces. Ruelle \cite{Ruelle1}
established the corresponding version in Hilbert spaces and
Ma\~n\'e~\cite{Man6} considered transformations in Banach spaces
under some compactness and invertibility assumptions, including
the case of differentiable maps with compact derivative at each
point. The results of Ma\~n\'e were extended by Thieullen in
\cite{Thie} for a class transformations satisfying a certain
asymptotic compactness. We refer the reader to the book by Hale,
Magalh\~aes and Oliva \cite{HMO} for a detailed discussion of the
state-of-the-art of the geometric theory of dynamical systems in
infinite dimensional spaces.

\subsection{Dynamical systems with nonzero Lyapunov exponents}\label{sec2.2}

The concept of hyperbolicity is closely related to the study of
Lyapunov exponents. These numbers measure the asymptotic
exponential rates of contraction and expansion in the neighborhood
of each given trajectory.

Let $f\colon M\to M$ be a diffeomorphism.  Given $x\in M$ and
$v\in T_xM$, we define the \emph{(forward) Lyapunov exponent} of
$(x, v)$ by
\[
\chi(x, v)=\limsup_{n\to+\infty}\frac1n\log\lVert d_xf^nv\rVert,
\]
with the convention that $\log0=-\infty$.  The abstract theory of
Lyapunov exponents (see \cite{lyap}), guarantees that for each
$x\in M$ there exist a positive integer $s(x)\le\dim M$, numbers
\[
\chi_1(x)<\cdots<\chi_{s(x)}(x),
\]
and linear spaces
\[
\{0\}=E_0(x)\subset E_1(x)\subset\cdots\subset E_{s(x)}(x)=T_xM
\]
such that if $i=1$, $\ldots$, $s(x)$ then
\[
E_i(x)=\{v\in T_xM:\chi(x,v)\le\chi_i(x)\},
\]
and $\chi(x,v)=\chi_i(x)$ whenever $v\in E_i(x)\setminus
E_{i-1}(x)$. Considering negative time, we can also define for
each $x\in M$ and $v\in T_xM$ the \emph{(backward) Lyapunov
exponent} of $(x,v)$ by
\[
\chi^-(x, v)=\limsup_{n\to-\infty}\frac1{\lvert n\rvert}\log\lVert
d_xf^nv\rVert.
\]
Again, the abstract theory of Lyapunov exponents guarantees that
for each $x\in M$ there exist a positive integer $s^-(x)\le\dim
M$, numbers
\[
\chi_1^-(x)>\cdots>\chi_{s^-(x)}^-(x), \]
and linear spaces
\[
T_xM=E_1^-(x)\supset\cdots\supset E_{s^-(x)}^-(x)\supset
E_{s^-(x)+1}^-(x)=\{0\}
\]
such that if $i=1$, $\ldots$, $s^-(x)$ then
\[
E_i^-(x)=\{v\in T_xM:\chi^-(x,v)\le\chi_i^-(x)\},
\]
and $\chi^-(x,v)=\chi_i^-(x)$ whenever $v\in E_i^-(x)\setminus
E_{i+1}^-(x)$. A~priori these two structures (for positive and
negative time) could be totally unrelated. The following result of
Oseledets \cite{Oseledec} shows that the two structures are indeed
related, in a very strong manner, in sets of full measure with
respect to any finite invariant measure.

\begin{theorem}[Multiplicative ergodic theorem]\label{Oset}
Let $f\colon M\to M$ be a $C^1$ diffeomorphism and $\mu$ an
$f$-invariant finite measure on~$M$ such that $\log^+\lVert
df\rVert$ and $\log^+\lVert df^{-1}\rVert$ are $\mu$-integrable.
Then for $\mu$-almost every point $x\in\Lambda$ there exist
subspaces $H_j(x)\subset T_xM$ for $j=1$, $\ldots$, $s(x)$ such
that:
\begin{enumerate}
\item
if $i=1$, $\ldots$, $s(x)$ then $E_i(x)=\bigoplus_{j=1}^iH_j(x)$ and
\[
\lim_{n\to\pm\infty}\frac1n\log\lVert d_xf^nv\rVert=\chi_i(x)
\]
with uniform convergence for $v$ on $\{v\in H_i(x):\lVert
v\rVert=1\}$; \item if $i\ne j$ then
\[
\lim_{n\to\pm\infty}\frac1n\log\lvert\angle(H_i(f^nx),H_j(f^nx))\rvert=0.
\]
\end{enumerate}
\end{theorem}

We note that if $M$ is a compact manifold then the functions
$\log^+\lVert df\rVert$ and $\log^+\lVert df^{-1}\rVert$ are
$\mu$-integrable for any finite measure $\mu$ on~$M$.  The
statement in Theorem~\ref{Oset} also holds in the more general
case of cocyles over a measurable transformation. See \cite{lyap}
for a detailed exposition and for a proof of the Multiplicative
ergodic theorem.

Let $f\colon M\to M$ be a $C^1$ diffeomorphism on a compact
manifold and $\mu$ an $f$-invariant finite Borel measure on~$M$.
We say that $f$ is \emph{nonuniformly hyperbolic} with respect to
$\mu$ if the set $\Lambda\subset M$ of points whose trajectories
are nonuniformly hyperbolic has measure $\mu(\Lambda)>0$.  In this
case the constants $\lambda$ and $\varepsilon$ in the definition
of nonuniformly hyperbolic trajectory are replaced by measurable
functions $\lambda(x)$ and $\varepsilon(x)$.

It follows from Theorem~\ref{Oset} that the following conditions
are equivalent:
\begin{enumerate}
\item $f$ is nonuniformly hyperbolic with respect to the measure
$\mu$; \item $\chi(x,v)\neq0$ for each $v\in T_xM$ and each $x$ in
a set with $\mu$-positive measure.
\end{enumerate}
Therefore, the nonuniformly hyperbolic diffeomorphisms with
respect to a given measure are precisely the diffeomorphisms with
all Lyapunov exponents nonzero in a set of positive measure.

One of the standing problems of the theory of nonuniformly
hyperbolic dynamical systems is to understand how ``common'' this
class is. Let $M$ be a compact smooth Riemannian manifold. It was
established by Katok in~\cite{Katok1} when $\dim M=2$ and by
Dolgopyat and Pesin in~\cite{DP} when $\dim M\ge3$ that there
exists a $C^\infty$ diffeomorphism $f$ such that:
\begin{enumerate}
\item $f$ preserves the Riemannian volume $m$ on $M$; \item $f$
has nonzero Lyapunov exponents at $m$-almost every point $x\in M$;
\item $f$ is a Bernoulli diffeomorphism.
\end{enumerate}
For any compact smooth Riemannian manifold $M$ of dimension at
least $5$, Brin constructed in \cite{Brin2} a $C^\infty$ Bernoulli
diffeomorphism which preserves the Riemannian volume and has all
but one Lyapunov exponent nonzero. On the other hand, the
construction of the above diffeomorphisms is not robust. In
another direction, Bochi \cite{Bochi} showed that on any compact
surface there exists a residual set $D$ of the $C^1$ area
preserving diffeomorphisms such each $f\in D$ is either an Anosov
diffeomorphism or has all Lyapunov exponents zero. This result was
announced by Ma\~n\'e but his proof was never published.

\section{Hyperbolic measures and dimension theory}\label{sec3}

\subsection{Product structure}\label{secff}

Let $f\colon M\to M$ be a diffeomorphism. We say that an
$f$-invariant measure $\mu$ on $M$ is a \emph{hyperbolic measure}
(with respect to~$f$) if all Lyapunov exponents are nonzero
$\mu$-almost everywhere, i.e., if $\chi(x,v)\neq0$ for each $v\in
T_xM$ and each $x$ in a set with full $\mu$-measure. One
consequence of the discussion in Section~\ref{sec2.2} is that the
existence of an invariant hyperbolic measure guarantees the
presence of nonuniform hyperbolicity and thus of a considerable
structure (see also Section~\ref{sec2.1}).  In this section we
describe in detail the structure of the hyperbolic measures.

Let $\mu$ be an $f$-invariant hyperbolic measure.  By
Theorem~\ref{thm:var}, for $\mu$-almost every point $x\in M$ there
exist local stable and unstable manifolds $V^s(x)$ and $V^u(x)$.
These manifolds somehow reproduce the product structure present in
the case of diffeomorphisms with hyperbolic sets, but a priori it
is unclear whether a given hyperbolic measure imitates or not the
product structure. This problem became known as the
Eckmann--Ruelle conjecture, claiming that locally a hyperbolic
measure indeed imitates the product structure defined by stable
and unstable manifolds. Even though Eckmann and Ruelle apparently
never formulated the conjecture, their work \cite{ER} discusses
several related problems and played a fundamental role in the
development of the theory and as such the expression seems
appropriate.

In order to formulate a rigorous result related to the resolution
of the conjecture we need the families of conditional measures
$\mu^s_x$ and $\mu^u_x$ generated by certain measurable partitions
constructed by Ledrappier and Young in \cite{LY}, based on former
work of Ledrappier and Strelcyn in~\cite{LedrappierStrelcyn}.  As
shown by Rohklin, any measurable partition $\xi$ of $M$ has
associated a family of conditional measures \cite{Ro}: for
$\mu$-almost every point $x\in M$ there exists a probability
measure $\mu_x$ defined on the element $\xi(x)$ of $\xi$
containing~$x$.  Furthermore, the conditional measures are
characterized  by the following property: if $\BB_\xi$ is a
$\sigma$-subalgebra of the Borel $\sigma$-algebra generated by the
unions of elements of $\xi$ then for each Borel set $A\subset M$,
the function $x\mapsto\mu_x(A\cap \xi(x))$ is $\BB_\xi$-measurable
and
\[
\mu(A)=\int_A\mu_x(A\cap\xi(x))\,d\mu.
\]
In \cite{LY}, Ledrappier and Young obtained two measurable
partitions $\xi^s$ and $\xi^u$ of~$M$ such that for $\mu$-almost
every point $x\in M$ we have:
\begin{enumerate}
\item $\xi^s(x)\subset V^s(x)$ and $\xi^u(x)\subset V^u(x)$; \item
for some $\gamma=\gamma(x)>0$,
\[
\xi^s(x)\supset V^s(x)\cap
B(x,\gamma)\quad\text{and}\quad
\xi^u(x)\supset V^s(x)\cap
B(x,\gamma).
\]
\end{enumerate}
We denote by $\mu_x^s$ and $\mu_x^u$ the conditional measures
associated respectively with the partitions $\xi^s$ and $\xi^u$.
We represent by $B^s(x,r)\subset V^s(x)$ and $B^u(x,r)\subset
V^u(x)$ the open balls centered at~$x$ of radius $r$ with respect
to the distances induced respectively in $V^s(x)$ and $V^u(x)$.

The following result of Barreira, Pesin and Schmeling \cite{annals}
establishes in the affirmative the Eckmann--Ruelle conjecture.

\begin{theorem}[Product structure of hyperbolic measures]\label{thm:annals}
Let $f\colon M\to M$ be a $C^{1+\alpha}$ diffeomorphism, for some
$\alpha>0$, and $\mu$ an $f$-invariant finite measure on~$M$ with
compact support. If $\mu$ is hyperbolic then given $\delta>0$
there exists a set $\Lambda\subset M$ with
$\mu(\Lambda)>\mu(M)-\delta$ such that for each $x\in\Lambda$ we
have
\[
r^\delta\le \frac{\mu(B(x,r))}{\mu^s_x(B^s(x,r))\mu^u_x(B^u(x,r))}\le
r^{-\delta}
\]
for all sufficiently small $r>0$.
\end{theorem}

This result was previously unknown even in the case of Anosov
diffeomorphisms. One of the major difficulties in the approach to
the problem has to do with the regularity of the stable and
unstable foliations that in general are not Lipschitz. In fact,
Schmeling showed in~\cite{Joerg2} that for a generic
diffeomorphism with a hyperbolic set, in some open set of
diffeomorphisms, the stable and unstable foliations are only
H\"older. Furthermore, the hyperbolic measure may not possess a
``uniform'' product structure even if the support does.

We know that the hypotheses in Theorem~\ref{thm:annals} are in a
certain sense optimal: Ledrappier and Misiurewicz \cite{LMi}
showed that the hyperbolicity of the measure is essential, while
Pesin and Weiss \cite{PW2} showed that the statement in
Theorem~\ref{thm:annals} cannot be extended to H\"older
homeomorphisms. On the other hand one does not know what happens
for $C^1$ diffeomorphisms that are not of class $C^{1+\alpha}$ for
some $\alpha>0$, particularly due to the nonexistence of an
appropriate theory of nonuniformly hyperbolic dynamical systems of
class~$C^1$.

\subsection{Dimension theory}\label{sec:dd}

There is a very close relation between the results described in
Section~\ref{secff} and the dimension theory of dynamical systems.
In order to describe this relation we briefly introduce  some
basic notions of dimension theory.

Let $X$ be a separable metric space. Given a set $Z\subset X$ and a
number $\alpha\in\RR$ we define
\[
m(Z,\alpha)=\lim_{\varepsilon\to 0} \inf_{\UU}
\sum_{U\in\UU}(\diam U)^\alpha,
\]
where the infimum is taken over all finite or countable covers of
$Z$ composed of open sets with diameter at most $\varepsilon$.
The \emph{Hausdorff dimension} of $Z$ is defined by
\[
\dim_HZ=\inf\{\alpha :m(Z,\alpha)=0\}.
\]
The \emph{lower} and \emph{upper box dimensions} of $Z$ are defined
respectively by
\[
\underline\dim_B Z= \liminf_{\varepsilon\to 0} \frac{\log
N(Z,\varepsilon)}{-\log\varepsilon} \quad\text{and}\quad
\overline\dim_B Z= \limsup_{\varepsilon\to 0} \frac{\log
N(Z,\varepsilon)}{-\log\varepsilon},
\]
where $N(Z,\varepsilon)$ denotes the number of balls of radius
$\varepsilon$ needed to cover $Z$.  It is easy to show that
\begin{equation}\label{eq:Z}
\dim_HZ\leq\underline\dim_B Z\leq\overline\dim_B Z.
\end{equation}
In general these inequalities can be strict and the coincidence of
the Hausdorff dimension and of the lower and upper box dimensions
is a relatively rare phenomenon that occurs only in some ``rigid''
situations (see Sections~\ref{sec4} and~\ref{sec5}; see \cite{p,
etds} for more details).

Let now $\mu$ be a finite measure on~$X$.  The \emph{Hausdorff
dimension} and the \emph{lower} and \emph{upper box dimensions} of
$\mu$ are defined respectively by
\[
\dim_H\mu= \lim_{\delta\to 0}\inf \{\dim_H Z:\mu(Z)\geq
\mu(X)-\delta\},
\]
\[
\underline\dim_B\mu= \lim_{\delta\to 0}\inf \{\underline\dim_B
Z:\mu(Z)\geq \mu(X)-\delta\},
\]
\[
\overline\dim_B\mu=\lim_{\delta\to 0}\inf \{\overline\dim_B
Z:\mu(Z)\geq \mu(X)-\delta\}.
\]
In general these quantities do not coincide, respectively, with
the {Hausdorff dimension} and the {lower} and {upper box
dimensions} of the support of $\mu$, and thus contain additional
information about the way in which the measure $\mu$ is
distributed on its support.  It follows immediately from
\eqref{eq:Z} that
\begin{equation}\label{eq:Z2}
{\dim_H}\mu\leq\underline\dim_B\mu\leq\overline\dim_B\mu.
\end{equation}
As with the inequalities in \eqref{eq:Z}, the inequalities in
\eqref{eq:Z2} are also strict in general. The following criterion
for equality was established by Young in \cite{Y}: if $\mu$ is a
finite measure on~$X$ and
\begin{equation}\label{pdim}
\lim_{r\to0}\frac{\log\mu(B(x,r))}{\log r}=d
\end{equation}
for $\mu$-almost every $x\in X$ then
\[
\dim_H\mu=\underline\dim_B\mu=\overline\dim_B \mu=d.
\]
The limit in \eqref{pdim}, when it exists, is called
\emph{pointwise dimension} of $\mu$ at~$x$.

In order to simplify the exposition we will assume that $\mu$ is
an ergodic measure, i.e., that any set $A\subset M$ such that
$f^{-1}A=A$ satisfies $\mu(A)=0$ or $\mu(M\setminus A)=0$. There
is in fact no loss of generality (see for example \cite{LY,
annals} for details). The following was established by Ledrappier
and Young in \cite{LY}.

\begin{theorem}[Existence of pointwise dimensions]\label{thm:LY}
Let $f$ be a $C^2$ diffeomorphism and $\mu$ an ergodic
$f$-invariant finite measure with compact support. If~$\mu$ is
hyperbolic then there exist constants $d^s$ and $d^u$ such that
\begin{equation}\label{dsds}
\lim_{r\to0}\frac{\log\mu_x^s(B^s(x,r))}{\log r}=d^s
\quad\text{and}\quad \lim_{r\to0}\frac{\log\mu_x^u(B^u(x,r))}{\log
r}=d^u
\end{equation}
for $\mu$-almost every point $x\in M$.
\end{theorem}

The limits in \eqref{dsds}, when they exist, are called
respectively \emph{stable} and \emph{unstable pointwise
dimensions} of $\mu$ at~$x$.

It was also established in \cite{LY} that
\begin{equation}\label{LY:in}
\limsup_{r\to0}\frac{\log\mu(B(x,r))}{\log r}\le d^s+d^u
\end{equation}
for $\mu$-almost every $x\in M$. It should be noted that
Ledrappier and Young consider a more general class of measures
in~\cite{LY}, for which some Lyapunov exponents may be zero. On
the other hand, they require the diffeomorphism $f$ to be of class
$C^2$. The only place in \cite{LY} where $f$ is required to be of
class $C^2$ concerns the Lipschitz regularity of the holonomies
generated by the intermediate foliations (such as any strongly
stable foliation inside the stable one). In the case of hyperbolic
measures a new argument was given by Barreira, Pesin and Schmeling
in \cite{annals} establishing the Lipschitz regularity for
$C^{1+\alpha}$ diffeomorphisms. This ensures that \eqref{dsds} and
\eqref{LY:in} hold almost everywhere even when $f$ is only of
class $C^{1+\alpha}$. See \cite{annals} for details.

For an ergodic hyperbolic finite measure with compact support that
is invariant under a $C^{1+\alpha}$ diffeomorphism,
Theorems~\ref{thm:annals} and~\ref{thm:LY} (and the above
discussion) imply that
\begin{equation}\label{op}
\lim_{r\to0}\frac{\log\mu(B(x,r))}{\log r}=d^s+d^u
\end{equation}
for $\mu$-almost every $x\in M$. Therefore, the above criterion by
Young allows us to conclude that
\[
{\dim_H}\mu=\underline\dim_B\mu=\overline\dim_B\mu=d^s+d^u.
\]
In fact the almost everywhere existence of the limit in \eqref{op}
guarantees the coincidence not only of these three dimensional
characteristics but also of many other characteristics of
dimensional type (see \cite{Y, p, annals} for more details).  This
allows us to choose any of these dimensional characteristics
according to the convenience in each application, since the common
value is always $d^s+d^u$.  Therefore, for hyperbolic measures
Theorems~\ref{thm:annals} and~\ref{thm:LY} allow a rigorous
approach to a ``fractal'' dimension of invariant measures, that is
well adapted to applications. Furthermore, the almost everywhere
existence of the limit in \eqref{op} plays the corresponding role
in dimension theory to the role of the Shannon--McMillan--Breiman
theorem in the entropy theory (see Section~\ref{sec:general}).

The $\mu$-almost everywhere existence of the limit in \eqref{op}
was established by Young \cite{Y} when $M$ is a surface and by
Ledrappier \cite{L} when $\mu$ is an SRB-measure (after Sinai,
Ruelle and Bowen; see for example \cite{lyap} for the definition).
In~\cite{PY}, Pesin and Yue extended the approach of Ledrappier to
hyperbolic measures with a ``quasi-product'' structure.
Theorem~\ref{thm:annals} shows that any hyperbolic measure
possesses a structure that is very close to the ``quasi-product''
structure. In \cite{str} (see also \cite{strextra}) Schmeling and
Troubetzkoy obtained versions of Theorem~\ref{thm:annals}
and~\eqref{op} for a class of endomorphisms.

\section{Dimension theory and thermodynamic formalism}\label{sec4}

\subsection{Dimension theory of geometric constructions}\label{seccc}

As we mentioned in Section~\ref{sec:dd} there are important
differences between the dimension theory of invariant sets and the
dimension theory of invariant measures. In~particular, while
virtually all dimensional characteristics of invariant hyperbolic
\emph{measures} coincide, the study of the dimension of invariant
hyperbolic \emph{sets} revealed that the different dimensional
characteristics frequently depend on other properties, and in
particular on number-theoretical properties. This justifies the
interest in simpler models in the context of the theory of
dynamical systems. We now make a little digression into the theory
of geometric constructions that precisely provides these models.

\begin{figure}[htbp]
\begin{psfrags}
\psfrag{1}{$\Delta_1$} \psfrag{2}{$\Delta_2$}
\psfrag{11}{$\Delta_{11}$} \psfrag{12}{$\Delta_{12}$}
\psfrag{21}{$\Delta_{21}$} \psfrag{22}{$\Delta_{22}$}
\begin{center}
\includegraphics{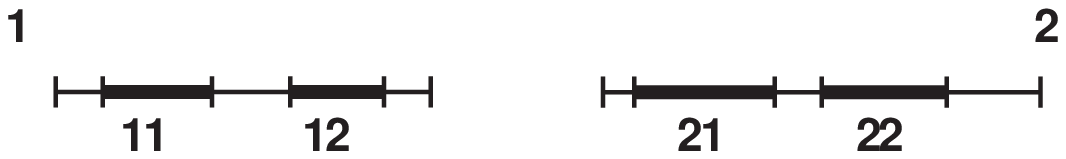}
\end{center}
\end{psfrags}
\caption{Geometric construction in $\RR$}\label{fig3}
\end{figure}

We start with the description of a geometric construction in
$\RR$. We consider constants $\lambda_1$, $\ldots$,
$\lambda_p\in(0,1)$ and disjoint closed intervals $\Delta_1$,
$\ldots$, $\Delta_p\subset\RR$ with length $\lambda_1$, $\ldots$,
$\lambda_p$ (see Figure~\ref{fig3}). For each $k=1$, $\ldots$,
$p$, we choose again $p$ disjoint closed intervals $\Delta_{k1}$,
$\ldots$, $\Delta_{kp}\subset\Delta_k$ with length ${\lambda_k
\lambda_1}$, $\ldots$, ${\lambda_k \lambda_p}$.  Iterating this
procedure, for each $n\in\NN$ we obtain $p^n$ disjoint closed
intervals $\Delta_{i_1\cdots i_n}$ with length $\prod_{k=1}^n
\lambda_{i_k}$. We define the set
\begin{equation}\label{ap4:eq1}
F=\bigcap_{n=1}^\infty\bigcup_{i_1\cdots i_n}\Delta_{i_1\cdots i_n}.
\end{equation}
In \cite{Mo}, Moran showed that $\dim_HF=s$ where $s$ is the unique
real number satisfying the identity
\begin{equation}\label{ap4:eq2}
\sum_{k=1}^p{\lambda_k}^s=1.
\end{equation}
It is remarkable that the Hausdorff dimension of $F$ does not
depend on the location of the intervals $\Delta_{i_1\cdots i_n}$
but only on their length.  Pesin and Weiss \cite{PW2} extended the
result of Moran to arbitrary symbolic dynamics in~$\RR^m$, using
the thermodynamic formalism (see Section~\ref{sec:for}).

To model hyperbolic invariant sets, we need to consider geometric
constructions described in terms of symbolic dynamics.  Given an
integer $p>0$, we consider the family of sequences
$X_p=\{1,\ldots,p\}^\NN$ and equip this space with the distance
\begin{equation}\label{g2}
d(\omega,\omega')=\sum_{k=1}^\infty
e^{-k}\lvert\omega_k-\omega'_k\rvert.
\end{equation}
We consider the shift map $\sigma\colon X_p\to X_p$ such that
$(\sigma\omega)_n=\omega_{n+1}$ for each $n\in\NN$.
A~\emph{geometric construction} in $\RR^m$ is defined by:
\begin{enumerate}
\item
a compact set $Q\subset X_p$ such that $\sigma^{-1}Q\supset Q$ for
some $p\in\NN$;
\item
a decreasing sequence of compact sets $\Delta_{\omega_1\cdots
\omega_n}\subset\RR^m$ for each $\omega\in Q$ with diameter
$\diam\Delta_{\omega_1\cdots \omega_n}\to 0$ as $n\to\infty$.
\end{enumerate}
We also assume that
\[
\interior\Delta_{i_1\cdots i_n}\cap\interior\Delta_{j_1\cdots
j_n}\ne\varnothing
\]
whenever $(i_1\cdots i_n)\ne(j_1\cdots j_n)$.
We define the \emph{limit set} $F$ of the geometric construction
by \eqref{ap4:eq1} with the union taken over all vectors
$(i_1,\ldots,i_n)$ such that $i_k=\omega_k$ for each $k=1$,
$\ldots$, $n$ and some $\omega\in Q$.

The geometric constructions include as a particular case the
iterated functions systems, that have been one of the main objects
of study of dimension theory, unfortunately sometimes with
emphasis on the form and not on the content. The situation
considered here has in mind applications to the dimension theory
and the multifractal analysis of dynamical systems (see the
following sections for a detailed description).

\begin{figure}[htbp]
\begin{psfrags}
\psfrag{1}{$\Delta_1$} \psfrag{2}{$\Delta_2$}
\psfrag{11}{$\Delta_{11}$} \psfrag{12}{$\Delta_{12}$}
\psfrag{21}{$\Delta_{21}$} \psfrag{22}{$\Delta_{22}$}
\begin{center}
\includegraphics{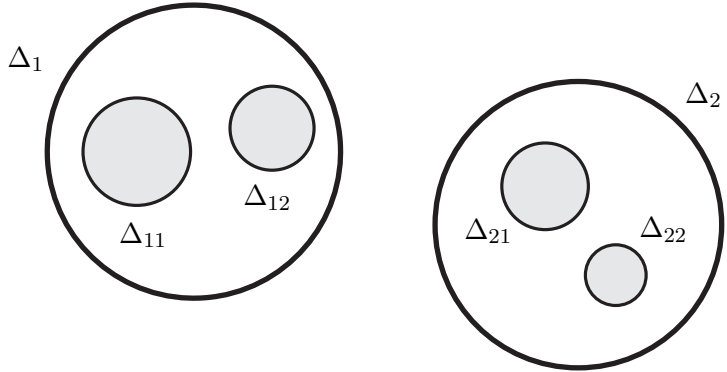}
\end{center}
\end{psfrags}
\caption{Geometric construction with balls}\label{fig4}
\end{figure}

We now consider the case in which all the sets $\Delta_{i_1\cdots
i_n}$ are balls (see Figure~\ref{fig4}). Write $r_{i_1\cdots
i_n}=\diam\Delta_{i_1\cdots i_n}$. The following result was
established by Barreira in~\cite{etds}.

\begin{theorem}[Dimension of the limit set]\label{thm:c}
For a geometric construction modelled by $Q\subset X_p$ for which
the sets $\Delta_{i_1\cdots i_n}$ are balls, if there exists a
constant $\delta>0$ such that
\[
r_{i_1\cdots i_{n+1}}\ge\delta r_{i_1\cdots
i_n}\quad\text{and}\quad r_{i_1\cdots i_{n+m}}\le r_{i_1\cdots
i_n}r_{i_{n+1}\cdots i_m}
\]
for each $(i_1i_2\cdots)\in Q$ and each $n$, $m\in\NN$ then
\[
\dim_H F=\underline\dim_B F=\overline\dim_B F=s,
\]
where $s$ is the unique real number satisfying the identity
\begin{equation}\label{ap4:eq3h}
\lim_{n\to\infty}\frac1n\log\sum_{i_1\cdots i_n}{ r_{i_1\cdots
i_n}}^s=0.
\end{equation}
\end{theorem}

We observe that this result contains as particular cases the
results of Moran and of Pesin and Weiss mentioned above (see also
Section~\ref{sec:for}), for which
\[
r_{i_1\cdots i_n}=\prod_{k=1}^n\lambda_{i_k}.
\]
The value of the dimension is also
independent of the location of the sets $\Delta_{i_1\cdots i_n}$.
We note that the hypotheses in Theorem~\ref{thm:c} naturally occur
in a class of invariant sets of uniformly hyperbolic dynamics (see
Section~\ref{bbb}).

\begin{figure}[htbp]
\begin{psfrags}
\psfrag{1}{$\Delta_1$} \psfrag{2}{$\Delta_2$}
\psfrag{11}{$\Delta_{11}$} \psfrag{12}{$\Delta_{12}$}
\psfrag{21}{$\Delta_{21}$} \psfrag{22}{$\Delta_{22}$}
\begin{center}
\includegraphics{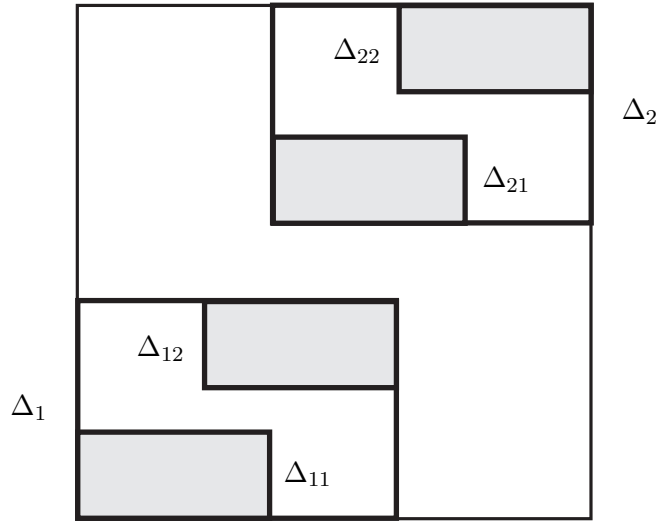}
\end{center}
\end{psfrags}
\caption{Number-theoretical properties and dimension
theory}\label{fig5}
\end{figure}

We now illustrate with an example how certain number-theoretical
properties can be relevant in dimension theory. We consider a
geometric construction in $\RR^2$ for which the sets
\[
\Delta_{i_1\cdots i_n}=(f_{i_1}\circ\cdots\circ
f_{i_n})([0,1]\times[0,1])
\]
are rectangles of sides $a^n$ and $b^n$, obtained through the
composition of the functions
\[
f_1(x,y)=(ax,by) \quad\text{and}\quad f_2(x,y)=(ax-a+1,by-b+1),
\]
where $a$, $b\in(0,1)$ with $b<1/2$ (see Figure~\ref{fig5}). In
particular, the projection of $\Delta_{i_1\cdots i_n}$ on the
horizontal axis is an interval with right endpoint given by
\begin{equation}\label{kkk4}
a^n+\sum_{k=0}^{n-1}j_ka^k,
\end{equation}
where $j_k=0$ if $i_k=1$ and $j_k=1-a$ if $i_k=2$. We assume now
that $a=(\sqrt5-1)/2$. In this case we have $a^2+a=1$ and thus for
each $n>2$ there exist several combinations $(i_1\cdots i_n)$ with
the same value in~\eqref{kkk4}. This duplication causes a larger
concentration of the sets $\Delta_{i_1\cdots i_n}$ in certain
regions of the limit set $F$. Thus, in view of computing the
Hausdorff dimension of $F$, when we take an open cover (see
Section~\ref{sec:dd}) it may be possible to replace, in the
regions of larger concentration of the sets $\Delta_{i_1\cdots
i_n}$, several elements of the cover by a unique element. This
procedure can cause $F$ to have a smaller Hausdorff dimension than
expected (with respect to the generic value obtained by Falconer
in \cite{FF1}). This was established by Neunh\"auserer
in~\cite{neuni}. See also \cite{PU, PoW} for former related
results. Additional complications can occur when $f_1$ and $f_2$
are replaced by functions that are not affine.

\subsection{Thermodynamic formalism}\label{sec:for}

The proof of Theorem~\ref{thm:c} is based on a ``nonadditive''
version of the topological pressure.  We first briefly introduce
the classical concept of topological pressure.

Given a compact set $Q\subset X_p$ such that $\sigma^{-1}Q\supset
Q$ and a continuous function $\varphi\colon Q\to\RR$ we define the
\emph{topological pressure} of $\varphi$ (with respect
to~$\sigma$) by
\begin{equation}\label{ap4:eq5}
P(\varphi)=\lim_{n\to\infty}\frac1n\log\sum_{i_1\cdots
i_n}\exp\sup\left(\sum_{k=0}^{n-1}\varphi\circ\sigma^k\right),
\end{equation}
where the supremum is taken over all sequences $(j_1 j_2\cdots)\in
Q$ such that $(j_1\cdots j_n)=(i_1\cdots i_n)$.  We define the
\emph{topological entropy} of $\sigma|Q$ by
\[
h(\sigma|Q)=P(0). \]
One can easily verify that
\[
h(\sigma|Q)=\lim_{n\to\infty}\frac1n\log N(Q,n),
\]
where $N(Q,n)$ is the number of vectors $(i_1,\ldots, i_n)$ whose
components constitute the first $n$ entries of some element of $Q$.

The topological pressure is one of the fundamental notions of the
thermodynamic formalism developed by Ruelle.  The topological
pressure was introduced by Ruelle in~\cite{R1} for expansive
transformations and by Walters in~\cite{W} in the general case.
For more details and references see \cite{bbook, KH, rbb,
walters}.

We now present an equivalent description of the topological
pressure. Let $\mu$ be a $\sigma$-invariant measure on $Q$ and
$\xi$ a countable partition of $Q$ into measurable sets. We write
\[
H_\mu(\xi)=-\sum_{C\in\xi}\mu(C)\log\mu(C),
\]
with the convention that $0\log0=0$.  We define the
\emph{Kolmogorov--Sinai entropy} of $\sigma|Q$ with respect to $\mu$
by
\begin{equation}\label{new:new}
h_\mu(\sigma|Q)=\sup_\xi\lim_{n\to\infty}\frac1n H_\mu\left(
\bigvee_{k=0}^{n-1}\sigma^{-k}\xi\right),
\end{equation}
where $\bigvee_{k=0}^{n-1}\sigma^{-k}\xi$ is the partition of $Q$
into sets of the form
\begin{equation}\label{tote}
C_{i_1\cdots i_n}=\bigcap_{k=0}^{n-1}\sigma^{-k}C_{i_{k+1}}
\end{equation}
with $C_{i_1}$, $\ldots$, $C_{i_n}\in\xi$ (it can be shown that
there exists the limit when $n\to\infty$ in~\eqref{new:new}). The
topological pressure satisfies the \emph{variational principle}
(see \cite{walters, KH} for details and references)
\begin{equation}\label{tytu}
P(\varphi)=\sup_\mu\left\{h_\mu(\sigma|Q)+\int_Q\varphi\,d\mu\right\},
\end{equation}
where the supremum is taken over all $\sigma$-invariant
probability measures on~$Q$.  A~$\sigma$-invariant probability
measure on $Q$ is called an \emph{equilibrium measure}
for~$\varphi$ (with respect to $\sigma|Q$) if the supremum in
\eqref{tytu} is attained by this measure, i.e., if
\[
P(\varphi)=h_\mu(\sigma|Q)+\int_Q\varphi\,d\mu.
\]

There exists a very close relation between dimension theory and
the thermodynamic formalism.  To illustrate this relation we
consider numbers $\lambda_1$,~$\ldots$, $\lambda_p$ and define the
function $\varphi\colon Q\to\RR$ by
\begin{equation}\label{pes}
\varphi(i_1 i_2\cdots)=\log \lambda_{i_1}.
\end{equation}
We have
\[
\begin{split}
P(s\varphi) &=\lim_{n\to\infty}\frac1n\log\sum_{i_1\cdots
i_n}\exp\left(s\sum_{k=1}^n\log\lambda_{i_k}\right)\\
&=\lim_{n\to\infty}\frac1n\log\sum_{i_1\cdots
i_n}\prod_{k=1}^n{\lambda_{i_k}}^s\\
&=\lim_{n\to\infty}\frac1n\log\left(\sum_{i=1}^p{\lambda_i}^s\right)^n\\
&=\log\sum_{i=1}^p{\lambda_i}^s.
\end{split}
\]
Therefore, the equation in \eqref{ap4:eq2} is equivalent to the
following equation involving the topological pressure:
\begin{equation}\label{bbeq}
P(s\varphi)=0.
\end{equation}
This equation was introduced by Bowen in~\cite{bqc} and is called
\emph{Bowen equation} (although it would be appropriate to call it
instead \emph{Bowen--Ruelle equation}). It has a rather universal
character: virtually all known equations to compute or estimate
the dimension of invariant sets of dynamical systems are
particular cases of this equation or of appropriate
generalizations. For example, the result of Pesin and Weiss in
\cite{PW2} mentioned in Section~\ref{seccc} can be formulated as
follows.

\begin{theorem}[Dimension of the limit set]\label{Pesu}
For a geometric construction modelled by $Q\subset X_p$ in which
the sets $\Delta_{i_1\cdots i_n}$ are balls of diameter
$\prod_{k=1}^n\lambda_{i_k}$, we have
\[
\dim_H F=\underline\dim_B F=\overline\dim_B F=s,
\]
where $s$ is the unique real number satisfying $P(s\varphi)=0$
with $\varphi$ as in \eqref{pes}.
\end{theorem}

However, the classical topological pressure is not adapted to all
geometric constructions.  Namely, comparing the equations in
\eqref{ap4:eq3h} and \eqref{ap4:eq5} it appears that it would be
appropriate to replace the sequence of functions
\begin{equation}\label{**4}
\varphi_n=\sum_{k=0}^{n-1}\varphi\circ\sigma^k
\end{equation}
in \eqref{ap4:eq5} by the new sequence
\[
\psi_n=s\log\diam\Delta_{i_1\cdots i_n}.
\]
We note that while the sequence $\varphi_n$ satisfies the identity
\[
\varphi_{n+m}=\varphi_n+\varphi_m\circ\sigma^n,
\]
the new sequence $\psi_n$ may not satisfy any additivity between
its terms.  Due to technical problems related with the existence
of the limit in \eqref{ap4:eq5} for sets that are not necessarily
compact, a different approach was used by Barreira in \cite{etds}
to introduce the nonadditive topological pressure. It is similar
to the introduction of the Hausdorff dimension, and uses the
theory of Carath\'eodory characteristics developed by Pesin (see
\cite{p} for references and full details).

Let $(X,\rho)$ be a compact metric space and $f\colon X\to X$ a
continuous transformation.  Given an open cover $\UU$ of~$X$ we
denote by $\WW_n(\UU)$ the collection of vectors
$\textbf{U}=(U_0,\ldots, U_n)$ of sets ${U_0}$, $\ldots$,
${U_n}\in\UU$ and write $m(\textbf{U})=n$. For each
$\textbf{U}\in\WW_n(\UU)$ we define the open set
\[
X(\textbf{U})=\bigcap_{k=0}^nf^{-k} U_k.
\]

We consider a sequence of functions $\Phi=\{\varphi_n\colon
X\to\RR\}_{n\in\NN}$. For each $n\in\NN$ we define
\[
\gamma_n(\Phi,\UU)=\sup\{|\varphi_n(x)-\varphi_n(y)|:\text{$x$,~$y\in
X(\textbf{U})$ for some $\textbf{U}\in\WW_n(\UU)$}\}
\]
and assume that
\begin{equation}\label{**5}
\limsup_{\diam\UU\to0}\limsup_{n\to\infty}\frac{\gamma_n(\Phi,\UU)}{n}=0.
\end{equation}
In the additive case (i.e., when $\Phi$ is composed by continuous
functions obtained as in \eqref{**4}) the condition \eqref{**5} is
always satisfied.

Given $\textbf{U}\in\WW_n(\UU)$, we write
\[
\varphi(\textbf{U})= \begin{cases} \sup_{X(\textbf{U})}\varphi_n
&\text{if $X(\textbf{U})\neq\varnothing$}\\
-\infty&\text{otherwise} \end{cases}.
\]
Given $Z\subset X$ and
$\alpha\in\RR$ we define
\[
M({Z,\alpha,\Phi,\UU})=\lim_{n\to\infty}
\inf_\Gamma\sum_{\textbf{U}\in\Gamma}\exp\left(-\alpha
m(\textbf{U})+\varphi(\textbf{U})\right),
\]
where the infimum is taken over all finite and infinite countable
families $\Gamma\subset\bigcup_{k\ge n}\WW_k(\UU)$ satisfying
$\bigcup_{\textbf{U}\in\Gamma}X(\textbf{U})\supset Z$. We define
\[
P_Z(\Phi,\UU)=\inf\{\alpha:M({Z,\alpha,\Phi,\UU})=0\}.
\]
The following properties were established in \cite{etds}:
\begin{enumerate}
\item
there exists the limit
\[
P_Z(\Phi)=\lim_{\diam\UU\to0}P_Z(\Phi,\UU);
\]
\item if there are constants $c_1$, $c_2>0$ such that $c_1
n\le\varphi_n\le c_2 n$ for each $n\in\NN$ and $h(f)<\infty$, then
there exists a unique number $s\in\RR$ such that $P_Z(s\Phi)=0$;
\item if there exists a continuous function $\psi\colon X\to\RR$
such that
\[
\text{$\varphi_{n+1}-\varphi_n\circ f\to\psi$ uniformly on~$X$}
\]
then
\[
P_X(\Phi)=\sup_\mu\left\{h_\mu(f)+\int_X\psi\,d\mu\right\},
\]
where the supremum is taken over all $f$-invariant probability
measures on~$X$.
\end{enumerate}

We call $P_Z(\Phi)$ the \emph{nonadditive topological pressure} of
$\Phi$ on the set~$Z$ (with respect to~$f$).  The nonadditive
topological pressure is a generalization of the classical
topological pressure and contains as a particular case the
subadditive version introduced by Falconer in \cite{f} under more
restrictive assumptions. In the additive case we recover the
notion of topological pressure introduced by Pesin and Pitskel'
in~\cite{PP}. The quantity $P_Z(0)$ coincides with the notion of
topological entropy for noncompact sets introduced in~\cite{PP},
and is equivalent to the notion introduced by Bowen in~\cite{Bo1}
(see~\cite{p}).

The equation $P_Z(s\Phi)=0$ is a nonadditive version of Bowen's
equation in~\eqref{bbeq}.  In particular, one can show that the
equation \eqref{ap4:eq3h} is equivalent to $P_Q(s\Phi)=0$, where
$\Phi$ is the sequence of functions $\varphi_n\colon Q\to \RR$
defined by
\begin{equation}\label{tg5}
\varphi_n(i_1i_2\cdots)=\log\diam\Delta_{i_1\cdots i_n}.
\end{equation}

\section{Hyperbolic sets and dimension theory}\label{sec5}

\subsection{Repellers and symbolic dynamics}\label{bbb}

As we observed in Section~\ref{seccc}, one of the main motivations
for the study of geometric constructions is the study of the
dimension theory of invariant sets of dynamical systems. This
approach can be effected with the use of Markov partitions.

We first consider the case of expanding maps. These constitute a
noninvertible version of the diffeomorphisms with hyperbolic sets.
Let $g\colon M\to M$ be a differentiable map of a compact
manifold. We consider a $g$-invariant set $J\subset M$, i.e., a
set such that $g^{-1}J=J$. We say that $J$ is a \emph{repeller}
of~$g$ and that $g$ is an \emph{expanding map} on~$J$ if there
exist constants $c>0$ and $\beta>1$ such that
\[
\lVert d_xg^nv\rVert\ge c\beta^n\lVert v\rVert
\]
for each $n\in\NN$, $x\in J$ and $v\in T_xM$.

Let now $J$ be a repeller of the differentiable map $g\colon M\to
M$. A~finite cover of~$J$ by nonempty closed sets $R_1$, $\ldots$,
$R_p$ is called a \emph{Markov partition} of $J$ if:
\begin{enumerate}
\item $\overline{\interior R_i}=R_i$ for each $i$; \item
$\interior R_i\cap \interior R_j=\varnothing$ whenever $i\ne j$;
\item $g R_i\supset R_j$ whenever $g(\interior R_i)\cap\interior
R_j\ne\varnothing$.
\end{enumerate}
The interior of each set $R_i$ is computed with respect to the
topology induced on~$J$. Any repeller has Markov partitions with
arbitrarily small diameter (see~\cite{R2}).

We can now use Markov partitions to model repellers by geometric
constructions. Let $J$ be a repeller of $g\colon M\to M$ and let
$R_1$, $\ldots$, $R_p$ be the elements of a Markov partition
of~$J$. We define a $p\times p$ matrix $A=(a_{ij})$ with entries
\[
a_{ij}=
\begin{cases}
1&\text{if $g(\interior R_i)\cap\interior R_j\ne\varnothing$}\\
0&\text{if $g(\interior R_i)\cap\interior R_j=\varnothing$}
\end{cases}.
\]
Consider the space of sequences $X_p=\{1,\ldots,p\}^\NN$ and the
shift map $\sigma\colon X_p\to X_p$ (see Section~\ref{seccc}). We
call \emph{topological Markov chain} with \emph{transition
matrix}~$A$ to the restriction of $\sigma$ to the set
\[
X_A=\{(i_1i_2\cdots)\in X_p:\text{$a_{i_ni_{n+1}}=1$ for every
$n\in\NN$}\}.
\]
We recall that a transformation $g$ is \emph{topologically mixing
on~$J$} if given open sets $U$ and $V$ with nonempty intersection
with $J$ there exists $n\in\NN$ such that $g^mU\cap V\cap
J\ne\varnothing$ for every $m>n$. If $g$ is topologically mixing
on $J$ then there exists $k\in\NN$ such that $A^k$ has only
positive entries.

It is easy to show that one can define a coding map $\chi\colon
X_A\to J$ by
\begin{equation}\label{coding}
\chi(i_1i_2\cdots)=\bigcap_{k=0}^\infty g^{-k}R_{i_{k+1}}.
\end{equation}
Furthermore, $\chi$ is surjective, satisfies
\begin{equation}\label{y2}
\chi\circ\sigma=g\circ \chi
\end{equation}
(i.e., the diagram in Figure~\ref{figx} is commutative), and is
H\"older continuous (with respect to the distance in $X_p$
introduced in~\eqref{g2}).

\begin{figure}[htbp]
\[
\begin{CD}
X_A @>\sigma>>X_A\\
@V{\chi}VV @VV{\chi}V\\
J @>g>> J
\end{CD}
\]
\caption{Symbolic coding of a repeller}\label{figx}
\end{figure}

Even though in general $\chi$ is not invertible (although $\card
\chi^{-1}x\le p^2$ for every~$x$), the identity in \eqref{y2}
allows us to see $\chi$ as a dictionary transferring the symbolic
dynamics $\sigma|X_A$ (and often the results on the symbolic
dynamics) to the dynamics of $g$ on~$J$. In particular, the
function $\chi$ allows us to see each repeller as a geometric
construction (see Section~\ref{seccc}) defined by the sets
\[
\Delta_{i_1\cdots i_n}=\bigcap_{k=0}^{n-1}g^{-k}R_{i_{k+1}}.
\]

We say that $g$ is \emph{conformal} on $J$ if $d_xg$ is a multiple
of an isometry for every $x\in J$.  When $J$ is a repeller of a
conformal transformation of class $C^{1+\alpha}$ one can show that
there is a constant $C>0$ such that
\begin{equation}\label{xaxa2}
C^{-1}\prod_{k=0}^{n-1}\exp\varphi(g^kx)\le\diam \Delta_{i_1\cdots
i_n}\le C\prod_{k=0}^{n-1}\exp\varphi(g^kx)
\end{equation}
for every $x\in \Delta_{i_1\cdots i_n}$, where the function
$\varphi\colon J\to\RR$ is defined by
\begin{equation}\label{toto}
\varphi(x)=-\log\lVert d_xg\rVert.
\end{equation}
We now use the topological pressure defined by \eqref{ap4:eq5}
with $Q=X_A$.

\begin{theorem}[Dimension of conformal repellers]\label{kk}
If $J$ is a repeller of a $C^{1+\alpha}$ transformation~$f$, for
some $\alpha>0$, such that $f$ is conformal on $J$ then
\[
\dim_HJ=\underline\dim_B J=\overline\dim_B J=s,
\]
where $s$ is the unique real number such that $P(s\varphi)=0$.
\end{theorem}

Ruelle established in \cite{R2} that $\dim_HJ=s$ (under the
additional assumption that $g$ is topologically mixing on $J$).
The coincidence between the Hausdorff dimension and the box
dimensions is due to Falconer \cite{FF}. The result in
Theorem~\ref{kk} was independently extended to expanding
transformations of class $C^1$ by Gatzouras and Peres in \cite{GP}
and by Barreira in \cite{etds} using different approaches. Under
the additional assumption that $g$ is of class $C^{1+\alpha}$ and
topologically mixing on~$J$, it was also shown by Ruelle in
\cite{R2} that if $\mu$ is the unique equilibrium measure of
$-s\varphi$ then
\begin{equation}\label{me:max}
\dim_HJ=\dim_H\mu.
\end{equation}
His proof consists in showing that $\mu$ is equivalent to the
$s$-dimensional Hausdorff measure on $J$ (in fact with
Radon--Nikodym derivative bounded and bounded away from zero).

Using \eqref{xaxa2} we find that
\[
P(s\varphi)=P_{X_A}(s\Phi)
\]
with $\varphi$ as in \eqref{toto} and being $\Phi$ the sequence of
functions defined by \eqref{tg5}. Furthermore, the conformality of
$g$ on~$J$ allows us to show that even though the sets
$\Delta_{i_1\cdots i_n}$ may not be balls they essentially behave
as if they were. In fact one would be able to reproduce with
little changes the proof of Theorem~\ref{thm:c} (and also the
proof of Theorem~\ref{Pesu}) to establish Theorem~\ref{kk}.
Nevertheless, there is  a technical difficulty related to the
possible noninvertibility of the coding map $\chi$
(see~\eqref{coding}). More generally, under the assumptions in
Theorem~\ref{kk} it was shown by Schmeling in \cite{joerg3} that
\[
\dim_H(\chi^{-1}B)=\dim_HB
\]
for \emph{any} subset $B\subset J$, provided that $X_A$ is given
the distance induced by the distance on $J$, so that
\[
\diam C_{i_1\cdots
i_n}=\prod_{k=0}^{n-1}\exp(\varphi\circ\chi)(\sigma^k\omega)=\lVert
d_{\sigma\omega} g^n\rVert^{-1}
\]
for each $\omega=(i_1i_2\cdots)\in X_A$ (see \eqref{tote}).

\subsection{Dimension theory in hyperbolic dynamics}\label{secsecd}

We now move to the study of the dimension of hyperbolic sets.

Let $\Lambda$ be a hyperbolic set for a diffeomorphism $f\colon
M\to M$.  We consider the functions $\varphi_s\colon\Lambda\to\RR$
and $\varphi_u\colon\Lambda\to\RR$ defined by
\[
\varphi_s(x)=\log\lVert d_xf|E^s(x)\rVert \quad\text{and}\quad
\varphi_u(x)=-\log\lVert d_xf|E^u(x)\rVert.
\]
Recall that $\Lambda$ is said to be \emph{locally maximal} if
there is an open neighborhood~$U$ of $\Lambda$ such that
\[
\Lambda=\bigcap_{n\in\ZZ}f^nU.
\]
The following result is a version of Theorem~\ref{kk} in the case
of hyperbolic sets.

\begin{theorem}[Dimension of hyperbolic sets on surfaces]\label{pop}
If $\Lambda$ is a locally maximal compact hyperbolic set of a
$C^1$ surface diffeomorphism, and $\dim E^s(x)=\dim E^u(x)=1$ for
every $x\in\Lambda$, then
\[
\dim_H\Lambda=\underline\dim_B\Lambda=\overline\dim_B\Lambda=t_s+t_u,
\]
where $t_s$ and $t_u$ are the unique real numbers such that
\[
P(t_s\varphi_s)=P(t_u\varphi_u)=0.
\]
\end{theorem}

It follows from work of McCluskey and Manning \cite{MM} that
$\dim_H\Lambda=t_s+t_u$. The coincidence between the Hausdorff
dimension and the lower and upper box dimensions is due to Takens
\cite{T} for $C^2$ diffeomorphisms (see also \cite{refPT}) and to
Palis and Viana \cite{PV} in the general case. Barreira
\cite{etds} and Pesin \cite{p} presented new proofs of
Theorem~\ref{pop} entirely based  on the thermodynamic formalism.

One can also ask whether there is an appropriate generalization of
property \eqref{me:max} in the present context, that is, whether
there exists an invariant measure $\mu$ supported on $\Lambda$
that satisfies $\dim_H\Lambda=\dim_H\mu$. However, the answer is
``almost always'' negative. More precisely, McCluskey and Manning
\cite{MM} showed that such a measure exists if and only if there
exists a continuous function $\psi\colon \Lambda\to\RR$ such that
\[
t_s\varphi_s- t_u\varphi=\psi\circ f -f
\]
on $\Lambda$. By Livschitz theorem (see for example \cite{KH}),
this happens if and only~if
\[
\lVert d_xf|E^s(x)\rVert^{t_s}\lVert d_xf|E^u(x)\rVert^{t_u}=1
\]
for every $x\in\Lambda$ and every $n\in\NN$ such that $f^nx=x$.
One can instead ask whether the supremum
\[
\delta(f) =\sup\{\dim_H\nu:\text{$\nu$ is an $f$-invariant measure
on $\Lambda$}\}
\]
is attained. Any invariant measure attaining this supremum would
attain the maximal complexity from the point of view of dimension
theory. The main difficulty of this problem is that the map
$\nu\mapsto\dim_H\nu$ is not upper semi-continuous: simply
consider the sequence $(\nu+(n-1)\delta)/n$ where $\dim_H\nu>0$
and $\delta$ is an atomic measure. It was shown by Barreira and
Wolf in \cite{BW1} (also using results in~\cite{BW2}) that the
supremum is indeed attained and by an ergodic measure, that is,
\[
\delta(f) =\max\{\dim_H\nu:\text{$\nu$ is an ergodic $f$-invariant
measure on $\Lambda$}\}.
\]
See \cite{wolf} for a precursor result in the special case of
polynomial automorphisms of $\CC^2$.

We now sketch the proof of Theorem~\ref{pop}. McCluskey and
Manning showed in~\cite{MM} that
\begin{equation}\label{hpopo}
\dim_H(V^s_\varepsilon(x)\cap\Lambda)=t_s \quad\text{and}\quad
\dim_H(V^u_\varepsilon(x)\cap\Lambda)=t_u
\end{equation}
for each $x\in\Lambda$.  Furthermore, Palis and Viana showed in
\cite{PV} that
\begin{equation}\label{hpopo2}
\dim_H(V^s_\varepsilon(x)\cap\Lambda)=
\underline\dim_B(V^s_\varepsilon(x)\cap\Lambda)=
\overline\dim_B(V^s_\varepsilon(x)\cap\Lambda)
\end{equation}
and
\begin{equation}\label{hpopo3}
\dim_H(V^u_\varepsilon(x)\cap\Lambda)=
\underline\dim_B(V^u_\varepsilon(x)\cap\Lambda)=
\overline\dim_B(V^u_\varepsilon(x)\cap\Lambda)
\end{equation}
for each $x\in\Lambda$.  Using these results, the completion of
the proof of Theorem~\ref{pop} depends in a crucial way on the
fact that the stable and unstable manifolds have dimension~$1$. In
fact, the product structure $[\cdot,\cdot]$ restricted to
$(V^s_\varepsilon(x)\cap\Lambda)\times(V^u_\varepsilon(x)\cap\Lambda)$
is locally a H\"older homeomorphism with H\"older inverse (and in
general is not more than H\"older, for generic diffeomorphisms in
a given open set, in view of work of Schmeling \cite{Joerg2}; see
also \cite{SRR}). However, when the stable and unstable manifolds
have dimension~$1$, the product structure $[\cdot,\cdot]$ is
locally a Lipschitz homeomorphism with Lipschitz inverse (see for
example \cite{KH}). This allows us to conclude that
\[
\dim_H[V^s_\varepsilon(x)\cap\Lambda,V^u_\varepsilon(x)\cap\Lambda]=
\dim_H((V^s_\varepsilon(x)\cap\Lambda)\times(V^u_\varepsilon(x)\cap\Lambda)),
\]
with corresponding identities between the lower and upper box
dimensions. Since the inequalities
\[
\dim_HA+\dim_HB\le\dim_H(A\times B) \quad\text{and}\quad
\overline\dim_B(A\times B)\le\overline\dim_B A+\overline\dim_B B
\]
are valid for any subsets $A$ and $B$ of $\RR^m$, it follows from
the identities in \eqref{hpopo}, \eqref{hpopo2} and~\eqref{hpopo3}
that
\begin{equation}\label{xoxa}
\begin{split}
\dim_H[V^s_\varepsilon(x)\cap\Lambda,V^u_\varepsilon(x)\cap\Lambda]
&=\underline\dim_B[V^s_\varepsilon(x)\cap\Lambda,V^u_\varepsilon(x)\cap\Lambda]\\
&=\overline\dim_B[V^s_\varepsilon(x)\cap\Lambda,V^u_\varepsilon(x)\cap\Lambda]\\
&= t_s+t_u.
\end{split}
\end{equation}
On the other hand, since $\Lambda$ is locally maximal one can
easily show that $[x,y]\in\Lambda$ for all sufficiently close $x$,
$y\in\Lambda$, or simply
that
\[
[V^s_\varepsilon(x)\cap\Lambda,V^u_\varepsilon(x)\cap\Lambda]\subset\Lambda
\]
for all $x\in\Lambda$ and all sufficiently small $\varepsilon$.
Choosing points $x_1$, $x_2$, $\ldots$ in $\Lambda$ such that
\[
\Lambda=\bigcup_{n\in\NN}[V^s_\varepsilon(x_n)
\cap\Lambda,V^u_\varepsilon(x_n)\cap\Lambda],
\]
Theorem~\ref{pop} follows now immediately from~\eqref{xoxa}.

We say that $f\colon M\to M$ is \emph{conformal} on $\Lambda$ when
$d_xf|E^s(x)$ and $d_xf|E^u(x)$ are multiples of isometries for
every point $x\in\Lambda$ (for example, if $M$ is a surface and
$\dim E^s(x)=\dim E^u(x)=1$ for every $x\in\Lambda$ then $f$ is
conformal on $\Lambda$). The proof of Theorem~\ref{pop} given by
Pesin in \cite{p} includes the case of conformal diffeomorphisms
on manifolds of arbitrary dimension (the statement can also be
obtained from results in \cite{etds}). In this situation the
product structure is still locally a Lipschitz homeomorphism with
Lipschitz inverse (see \cite{Boris, p} for details) and thus we
can use the same approach as above.

The study of the dimension of repellers and hyperbolic sets of
nonconformal transformations is not yet as developed. The main
difficulty has to do with the possibility of existence of distinct
Lyapunov exponents associated to directions that may change from
point to point.  There exist however some partial results, for
certain classes of repellers and hyperbolic sets, starting
essentially with the seminal work by Douady and Oesterl\'e
in~\cite{DO}. Namely, Falconer \cite{FF2} computed the Hausdorff
dimension of a class of nonconformal repellers (see also
\cite{FF1}), while Hu \cite{hu} computed the box dimension of a
class of nonconformal repellers that leave invariant a strong
unstable foliation. Related ideas were applied by Simon and
Solomyak in \cite{SS} to compute the Hausdorff dimension of a
class of hyperbolic sets in~$\RR^3$. Falconer also studied a class
of limit sets of geometric constructions obtained from a
composition of affine transformations that are not necessarily
conformal \cite{FF1}. In another direction, Bothe~\cite{Bothe} and
Simon \cite{Simon} (also using his methods in \cite{Simon0} for
noninvertible transformations) studied the dimension of solenoids
(see \cite{p, sw} for details). A solenoid is a hyperbolic set of
the form $\Lambda=\bigcap_{n=1}^\infty f^nT$, where
$T\subset\RR^3$ is diffeomorphic to a ``solid torus'' $S^1\times
D$ for some closed disk $D\subset\RR^2$ and $f\colon T\to T$ is a
diffeomorphism such that for each $x\in S^1$ the intersection
$f(T)\cap(\{x\}\times D)$ is a disjoint union of $p$ sets
homeomorphic to a closed disk.

In a similar way that in Section~\ref{bbb} the proof of the
identities in \eqref{hpopo} can be obtained with the use of Markov
partitions. We briefly recall the notion of Markov partition for a
hyperbolic set. A~nonempty closed set $R\subset\Lambda$ is called
a \emph{rectangle} if $\diam R<\delta$ (where $\delta$ is given by
the product structure; see Section~\ref{sec1.1}),
$\overline{\interior R}=R$, and $[x,y]\in R$ whenever $x$, $y\in
R$. A~finite cover of $\Lambda$ by rectangles $R_1$, $\ldots$,
$R_p$ is called a \emph{Markov partition} of $\Lambda$ if:
\begin{enumerate}
\item $\interior R_i\cap \interior R_j=\varnothing$ whenever $i\ne
j$; \item if $x\in f(\interior R_i)\cap\interior R_j$ then
\[
f^{-1}(V^u_\varepsilon(fx)\cap R_j)\subset V^u_\varepsilon(x)\cap
R_i \quad\text{and}\quad f(V^s_\varepsilon(x)\cap R_i)\subset
V^s_\varepsilon(fx)\cap R_j.
\]
\end{enumerate}
The interior of each set $R_i$ is computed with respect to the
topology induced on~$\Lambda$. Any hyperbolic set possesses Markov
partitions with arbitrarily small diameter (see \cite{bbook} for
references and full details).

Let $\partition=\{R_1, \ldots, R_p\}$ be a Markov partition of a
hyperbolic set. It is well known that
$\partial\partition=\bigcup_{i=1}^p\partial R_i$ has zero measure
with respect to any equilibrium measure. This is a simple
consequence of the fact that $\partial\partition$ is a closed set
with dense complement. On the other hand, it is also interesting
to estimate the measure of neighborhoods of $\partial\partition$.
This can be simpler when each element of $\partition$ has a
piecewise regular boundary (as in the case of hyperbolic
automorphisms of $\TT^2$), but it is well known that Markov
partitions may have a very complicated boundary. In particular, it
was discovered by Bowen \cite{boggg} that $\partial\partition$ is
not piecewise regular in the case of hyperbolic automorphisms
of~$\TT^3$.  It was shown by Barreira and Saussol in
\cite{product} that if $\mu$ is an equilibrium measure of a
H\"older continuous function then there exist constants $c>0$ and
$\nu>0$ such that
\[
\mu(\{ x\in \Lambda\colon d(x,\partial \partition)<\varepsilon\})
\le c \varepsilon^\nu
\]
for every $\varepsilon>0$. This provides a control of the measure
near $\partial\partition$. Furthermore, it is possible to consider
any $\nu>0$ such that
\[
\nu<\frac{P_\Lambda(\varphi)-P_I( \varphi)}{\log\max\{\lVert
d_xf\rVert:x\in\Lambda\}},
\]
where $\varphi$ in chosen in such a way that $\mu$ is an
equilibrium measure of $\varphi$, and
\[
I=\bigcup_{n\in\ZZ}f^n(\partial\partition)
\]
is the invariant hull of $\partial\partition$.

\section{Multifractal analysis}\label{sec6}

\subsection{Hyperbolic dynamics and Birkhoff averages}\label{sec6.1}

As we described in Section~\ref{sec:dd}, if $f\colon M\to M$ is a
$C^{1+\alpha}$ diffeomorphism and $\mu$ is an ergodic
$f$-invariant hyperbolic finite measure, then
\begin{equation}\label{pdim2}
\lim_{r\to0}\frac{\log\mu(B(x,r))}{\log r}=d^s+d^u
\end{equation}
for $\mu$-almost every point $x\in M$, where the numbers $d^s$ and
$d^u$ are as in \eqref{dsds}. Of course that this does not mean
that \emph{all} points necessarily satisfy \eqref{pdim2}.
Multifractal analysis precisely studies  the properties of the
level sets
\begin{equation}\label{koko}
\left\{x\in M:\lim_{r\to0}\frac{\log\mu(B(x,r))}{\log
r}=\alpha\right\}
\end{equation}
for each $\alpha\in\RR$. In this section we present the main
components of multifractal analysis and describe its relation with
the theory of dynamical systems.

Birkhoff's ergodic theorem---one of the basic but fundamental
results of ergodic theory---states that if $S\colon X\to X$ is a
measurable transformation preserving a finite measure $\mu$
on~$X$, then for each function $\varphi\in L^1(X,\mu)$ the limit
\[
\varphi_S(x)=\lim_{n\to\infty}\frac1n\sum_{k=0}^{n-1}\varphi(S^kx)
\]
exists for $\mu$-almost every point $x\in X$. Furthermore, if
$\mu$ is ergodic (see Section~\ref{sec:dd}) then
\begin{equation}\label{ttg}
\varphi_S(x)=\frac1{\mu(X)}\int_X\varphi\,d\mu
\end{equation}
for $\mu$-almost every $x\in X$. Again this does not mean that the
identity in \eqref{ttg} is valid for \emph{every} point $x\in X$
for which $\varphi_S(x)$ is well-defined.  For each $\alpha\in\RR$
we define the level set
\[
K_\alpha(\varphi)=\left\{x\in
X:\lim_{n\to\infty}\frac1n\sum_{k=0}^{n-1}\varphi(S^kx)=\alpha\right\},
\]
i.e., the set of points $x\in X$ such that $\varphi_S(x)$ is
well-defined and equal to~$\alpha$.  We also consider the set
\begin{equation}\label{liv}
K(\varphi)=\left\{x\in
X:\liminf_{n\to\infty}\frac1n\sum_{k=0}^{n-1}\varphi(S^kx)
<\limsup_{n\to\infty}\frac1n\sum_{k=0}^{n-1}\varphi(S^kx)\right\}.
\end{equation}
It is clear that
\begin{equation}\label{kkk}
X=K(\varphi)\cup\bigcup_{\alpha\in\RR}K_\alpha(\varphi).
\end{equation}
Furthermore, the sets in this union (possibly uncountable) are
pairwise disjoint.  We call the decomposition of~$X$ in
\eqref{kkk} a \emph{multifractal decomposition}.

One way to measure the complexity of the sets $K_\alpha(\varphi)$
is to compute their Hausdorff dimension.  We define a function
\[
\DD\colon\{\alpha\in\RR:K_\alpha(\varphi)\ne\varnothing\}\to\RR
\]
by
\[
\DD(\alpha)=\dim_HK_\alpha(\varphi).
\]
We also define the numbers
\[
\underline\alpha=\inf\left\{\int_X\varphi\,d\mu:\mu\in\MM\right\}
\quad\text{and}\quad
\overline\alpha=\sup\left\{\int_X\varphi\,d\mu:\mu\in\MM\right\},
\]
where $\MM$ represents the family of $S$-invariant probability
measures on~$X$.  It is easy to verify that
$K_\alpha(\varphi)=\varnothing$ whenever
$\alpha\not\in[\underline\alpha,\overline\alpha]$.  We also define
the function $T\colon\RR\to\RR$ by
\[
T(q)=P(q\varphi)-qP(\varphi)
\]
(where $P$ denotes the topological pressure). For topological
Markov chains (see Section~\ref{bbb}) the function $T$ is analytic
(see the book by Ruelle \cite{rbb}). Under the assumptions in
Theorem~\ref{batata} below there exists a unique equilibrium
measure $\nu_q$ of $q\varphi$ (see Section~\ref{sec:for}).

The following result shows that in the case of topological Markov
chains the set $K_\alpha(\varphi)$ is nonempty for any
$\alpha\in(\underline\alpha,\overline\alpha)$ and that the
function $\DD$ is analytic and strictly convex.

\begin{theorem}[Multifractal analysis of Birkhoff averages]\label{batata}
If $\sigma|X$ is a to\-pologically mixing topological Markov chain
and $\varphi\colon X\to\RR$ is a H\"older continuous function
then:
\begin{enumerate}
\item\label{tt2} $K_\alpha(\varphi)$ is dense in~$X$ for each
$\alpha\in(\underline\alpha,\overline\alpha)$; \item the function
$\DD\colon(\underline\alpha,\overline\alpha)\to\RR$ is analytic
and strictly convex; \item the function $\DD$ is the Legendre
transform of $T$, i.e.,
\[
\DD(-T'(q))=T(q)-qT'(q)
\]
for each $q\in\RR$; \item if $q\in\RR$ then
$\nu_q(K_{-T'(q)}(\varphi))=1$ and
\[
\lim_{r\to0}\frac{\log \nu_q(B(x,r))}{\log r}=T(q)-qT'(q)
\]
for $\nu_q$-almost every point $x\in K_{-T'(q)}(\varphi)$.
\end{enumerate}
\end{theorem}

\begin{figure}[htbp]
\begin{psfrags}
\psfrag{Q}{$q$} \psfrag{T}{$T(q)$}
\begin{center}
\includegraphics{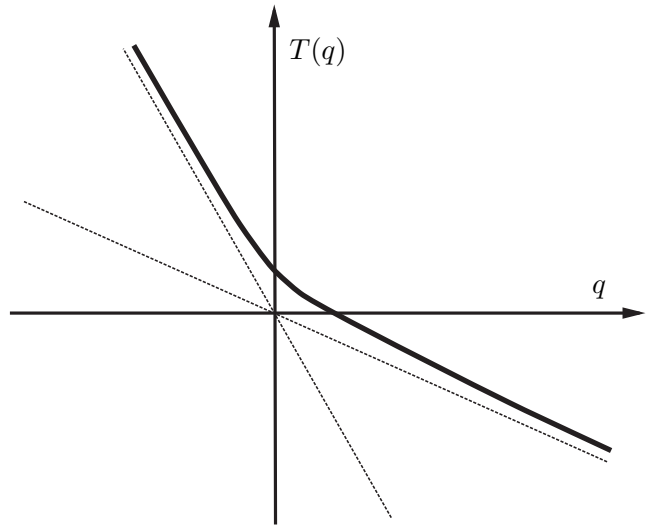}
\end{center}
\end{psfrags}
\caption{Graph of the function $T$}\label{fig6}
\end{figure}

\begin{figure}[htbp]
\begin{psfrags}
\psfrag{A}{$\alpha$} \psfrag{D}{$\DD(\alpha)$} \psfrag{Neg}{$q<0$}
\psfrag{BBB}{$q>0$}
\begin{center}
\includegraphics{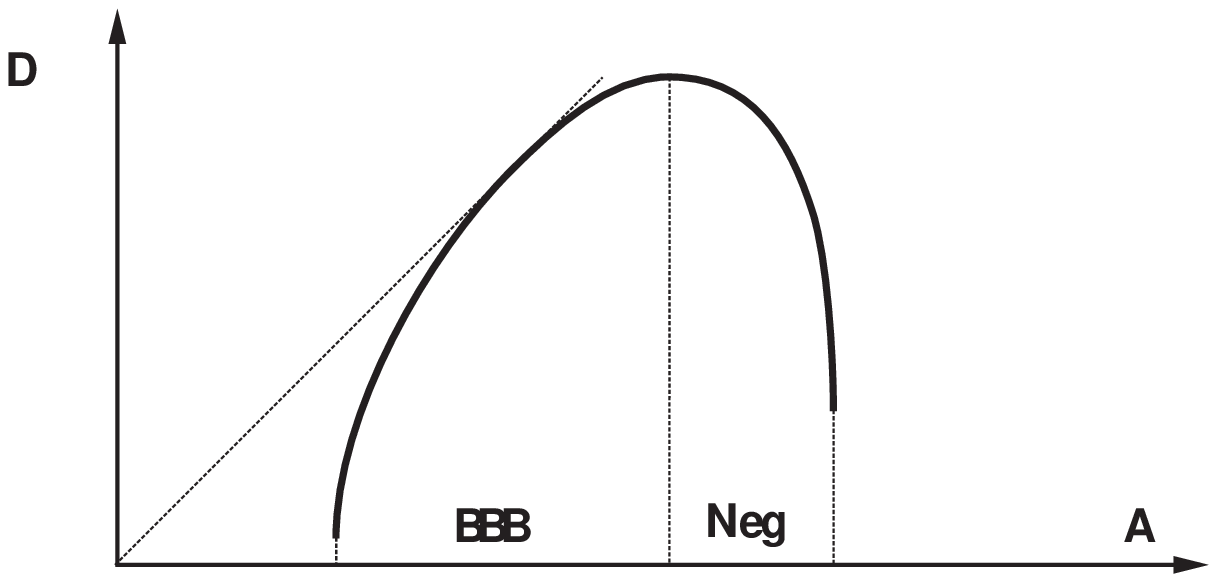}
\end{center}
\end{psfrags}
\caption{Graph of the function $\DD$}\label{fig7}
\end{figure}

See Figures~\ref{fig6} and~\ref{fig7} for typical graphs of the
functions $T$ and $\DD$.

Statement~\ref{tt2} in Theorem~\ref{batata} is an exercise (note
that we are considering one-sided topological Markov chains,
although all the results readily extend to the case of two-sided
topological Markov chains). The remaining statements in
Theorem~\ref{batata} are a particular case of results formulated
by Barreira and Schmeling in~\cite{ijm}. These were obtained as a
consequence of results of Pesin and Weiss in \cite{pw}, where they
effect a multifractal analysis for conformal repellers (see
Section~\ref{bbb}). In \cite{joerg}, Schmeling showed that the
domain of $\DD$ coincides with
$[\underline\alpha,\overline\alpha]$, i.e., that
$K_\alpha(\varphi)\ne\varnothing$ if and only if
$\alpha\in[\underline\alpha,\overline\alpha]$.

The concept of multifractal analysis was suggested by Halsey,
Jensen, Ka\-da\-noff, Procaccia and Shraiman in
\cite{HalseyJensenKadanoffProcacciaShraiman}. The first rigorous
approach is due to Collet, Lebowitz and Porzio in \cite{CLP} for a
class of measures invariant under one-dimensional Markov maps. In
\cite{lopes}, Lopes considered the measure of maximal entropy for
hyperbolic Julia sets, and in \cite{rand}, Rand studied Gibbs
measures for a class of repellers. We refer the reader to the book
by Pesin~\cite{p} for a more detailed discussion and further
references.

Theorem~\ref{batata} reveals an enormous complexity of
multifractal decompositions that is not foreseen by Birkhoff's
ergodic theorem. In particular it shows that the multifractal
decomposition in \eqref{kkk} is composed by an \emph{uncountable}
number of (pairwise disjoint) dense invariant sets, each of them
having positive Hausdorff dimension. We will see in
Section~\ref{sec7.1} that the set $K(\varphi)$ in \eqref{liv} is
also very complex (even though it has zero measure with respect to
\emph{any} finite invariant measure, as a simple consequence of
Birkhoff's ergodic theorem).

We now come back to the study of the level sets in \eqref{koko}.
Let $M$ be a surface and  $\Lambda\subset M$  a locally maximal
compact hyperbolic set for a $C^{1+\alpha}$ diffeomorphism
$f\colon M\to M$.  We assume that $f$ is topologically mixing on
$\Lambda$ (see Section~\ref{bbb}). Consider an equilibrium measure
$\mu$ of a H\"older continuous function
$\varphi\colon\Lambda\to\RR$. Under these assumptions $\mu$ is
unique and thus it is ergodic (see for example \cite{KH}).

We define functions $T_s\colon\Lambda\to\RR$ and
$T_u\colon\Lambda\to\RR$ by
\[
T_s(q)=P(-q\log\lVert df|E^s\rVert+q\varphi)-qP(\varphi)
\]
and
\[
T_u(q)=P(q\log\lVert df|E^u\rVert+q\varphi)-qP(\varphi).
\]
In \cite{simp}, Simpelaere showed that
\[
\dim_H\left\{x\in M:\lim_{r\to0}\frac{\log\mu(B(x,r))}{\log
r}=\alpha\right\} =T_s(q)-qT'_s(q)+T_u(q)-qT'_u(q),
\]
where $q\in\RR$ is the unique real number such that
$\alpha=-T'_s(q)-T'_u(q)$.  Another proof of this statement was
given by Pesin and Weiss in \cite{pwchaos} (see also \cite{p}).
Again we observe an enormous complexity that is not precluded by
the $\mu$-almost everywhere existence of the pointwise dimension.

In the case of hyperbolic flows versions of these results were
obtained by Barreira and Saussol \cite{flow} (in the case of
entropy spectra; see Section~\ref{sec:general}) and by Pesin and
Sadovskaya \cite{Psad}, using in particular the symbolic dynamics
developed by Bowen \cite{bttt} and Ratner \cite{ratner}.

\subsection{General concept of multifractal analysis}\label{sec:general}

In fact the approach of multifractal analysis extends to many
other classes of dynamical systems and to other local quantities.
With the purpose of unifying the theory, Barreira, Pesin and
Schmeling \cite{chaos} proposed a general concept of multifractal
analysis that we now describe.

We consider a function $g\colon Y\to[-\infty,+\infty]$ defined on
a subset $Y$ of~$X$.  The level sets
\[
K^g_\alpha=\{x\in X:g(x)=\alpha\}
\]
are pairwise disjoint and we obtain a \emph{multifractal
decomposition} of~$X$ given by
\begin{equation}\label{kkk2}
X=(X\setminus Y)\cup\bigcup_{\alpha\in[-\infty,+\infty]}K^g_\alpha.
\end{equation}
Let now $G$ be a function defined on the subsets of~$X$. We define
the \emph{multifractal spectrum}
$\FF\colon[-\infty,+\infty]\to\RR$ of the pair $(g,G)$ by
\[
\FF(\alpha)=G(K^g_\alpha).
\]

When $X$ is a compact manifold and $g$ is differentiable, each
level set $K_\alpha^g$ is a hyper-surface for all values of
$\alpha$ that are not critical values of $g$. In multifractal
analysis we are mostly interested in the study of level sets of
functions that are not differentiable (and typically are only
measurable), that  naturally appear in the theory of dynamical
systems.

Furthermore, the multifractal spectra encode precious information
about these functions and ultimately about the dynamical system
that originated them. In applications we have frequently no
information about the ``microscopic'' nature of the dynamical
system but only information about ``macroscopic'' quantities such
as for example about multifractal spectra. It is therefore
important to try to recover information about the dynamical system
through the information given by these ``macroscopic'' quantities
and in particular by the multifractal spectra (see also the
discussion in Section~\ref{sec8.2}).

We now describe some of the functions $g$ and $G$ that  naturally
occur in dynamical systems. Let $X$ be a separable metric space
and $f\colon X\to X$ a continuous function.  We define functions
$G_D$ and $G_E$ by
\[
G_D(Z)=\dim_HZ \quad\text{and}\quad G_E(Z)=h(f|Z).
\]
We call \emph{dimension spectra} and \emph{entropy spectra}
respectively to the spectra generated by $G_D$ and $G_E$.

Let $\mu$ be a finite Borel measure on~$X$ and $Y\subset X$ the
set of points $x\in X$ for which the pointwise dimension
\[
g_D(x)=g_D^{(\mu)}(x)=\lim_{r\to0}\frac{\log\mu(B(x,r))}{\log r}
\]
is well-defined. We obtain two multifractal spectra
\[
\DD_D=\DD_D^{(\mu)}\quad\text{and}\quad \DD_E=\DD_E^{(\mu)}
\]
specified respectively by the pairs $(g_D,G_D)$ and $(g_D,G_E)$.
For $C^{1+\alpha}$ diffeomorphisms and hyperbolic invariant
measures, Theorem~\ref{thm:annals} ensures that $\mu(X\setminus
Y)=0$.

Let now $X$ be a separable metric space and $f\colon X\to X$ a
continuous transformation preserving a probability measure $\mu$
on~$X$.  Given a partition $\xi$ of~$X$, for each $n\in\NN$ we
define a new partition of~$X$ by
$\xi_n=\bigvee_{k=0}^{n-1}f^{-k}\xi$.  We consider the set $Y$
formed by the points $x\in X$ for which the \emph{local entropy}
\[
g_E(x)=g_E^{(\mu)}(x)=\lim_{n\to\infty}-\frac1n\log\mu(\xi_n(x))
\]
is well-defined, where $\xi_n(x)$ denotes the element of $\xi_n$
containing~$x$. By the Shannon--McMillan--Breiman theorem in the
entropy theory (see for example~\cite{KH}) we have $\mu(X\setminus
Y)=0$. Furthermore, if $\xi$ is a generating partition and $\mu$
is ergodic then $g_E(x)=h_\mu(f)$ for $\mu$-almost every point
$x\in X$.  We obtain two multifractal spectra
\[
\EE_D=\EE_D^{(\mu)}\quad\text{and}\quad\EE_E=\EE_E^{(\mu)}
\]
specified respectively by the pairs $(g_E,G_D)$ and $(g_E,G_E)$.

We can also consider functions defined by the Lyapunov exponents.
In this case Theorem~\ref{Oset} guarantees that the involved
limits exist almost everywhere. We only consider a particular case
here. Let $X$ be a differentiable manifold and $f\colon X\to X$ a
$C^1$ map. Consider the set $Y\subset X$ of points $x\in X$ for
which the limit
\[
\lambda(x)=\lim_{n\to+\infty}\frac1n\log\lVert d_xf^n\rVert
\]
exists. By Theorem~\ref{Oset} (or by Kingman's subadditive ergodic
theorem), if $\mu$ is an $f$-invariant probability Borel measure,
then $\mu(X\setminus Y)=0$. We define the function $g_L$ on $Y$ by
\[
g_L(x)=\lambda(x).
\]
We obtain two multifractal spectra $\LL_D$ and $\LL_E$ specified
respectively by the pairs $(g_L,G_D)$ and $(g_L,G_E)$.

The spectrum $\DD_D=\DD$ was already considered in
Section~\ref{sec6.1}. We now describe the spectrum $\EE_E$.  Given
a compact hyperbolic set $\Lambda$ and a continuous function
$\varphi\colon\Lambda\to\RR$ we define the function
$T_E\colon\RR\to\RR$ by
\[
T_E(q)=P(q\varphi)-qP(\varphi).
\]
Under the assumptions in Theorem~\ref{EE} below there exists a
unique equilibrium measure $\nu_q$ of $q\varphi$. We also define
the numbers
\[
\underline\beta=\inf\left\{-\int_\Lambda\varphi\,d\mu:\mu\in\MM\right\}
\quad\text{and}\quad
\overline\beta=\sup\left\{-\int_\Lambda\varphi\,d\mu:\mu\in\MM\right\},
\]
where $\MM$ denotes the family of $f$-invariant probability
measures on $\Lambda$.

\begin{theorem}[Multifractal analysis of the spectrum $\EE_E$]\label{EE}
Let $\Lambda$ be a compact hyperbolic set of a $C^1$
diffeomorphism $f\colon M\to M$, such that $f$ is topologically
mixing on $\Lambda$. If $\mu$ is the equilibrium measure of a
H\"older continuous function $\varphi\colon\Lambda\to\RR$ then:
\begin{enumerate}
\item $K_\alpha^{g_E}$ is dense in $\Lambda$ for each
$\alpha\in(\underline\beta,\overline\beta)$; \item the function
$\EE_E$ is analytic and strictly convex on
$(\underline\beta,\overline\beta)$; \item the function $\EE_E$ is
the Legendre transform of $T_E$; \item for each $q\in\RR$ we have
\[
g_E^{(\nu_q)}(x)=T_E(q)-qT'_E(q)
\]
for $\nu_q$-almost every point
$x\in\Lambda$.
\end{enumerate}
\end{theorem}

Theorem~\ref{EE} is an immediate consequence of results of
Barreira, Pesin and Schmeling in \cite{jdcs} (see also
\cite{ijm}).

We note that in Theorem~\ref{EE} the manifold is not necessarily
two-dimen\-sional, contrarily to what happens in
Section~\ref{sec6.1} when we described the spectrum $\DD_D$. In
the case of conformal repellers, Pesin and Weiss \cite{pw}
obtained a multifractal analysis of $\DD_D$ and Barreira, Pesin
and Schmeling~\cite{chaos} obtained a multifractal analysis of
$\EE_E$. In \cite{tv2}, Takens and Verbitski obtained a
multifractal analysis of the spectrum $\EE_E$ for expansive
homeomorphisms with specification and a certain class of
continuous functions (note that these systems need not have Markov
partitions).

We note that the spectra $\DD_D$ and $\EE_E$ are of different
nature from that of the spectra $\DD_E$ and $\EE_D$. Namely, the
first two relate pointwise quantities---the pointwise dimension
and the local entropy---with global quantities that are naturally
associated to them---respectively the Hausdorff dimension and the
Kolmogorov--Sinai entropy.  On the other hand, the spectra $\DD_E
$ and $\EE_D$ mix local and global quantities of distinct nature.
We refer to them as \emph{mixed spectra}. It is reasonable to
expect that the mixed spectra contain additional information about
the dynamical system. It is also possible to describe the
multifractal properties of these spectra although this requires a
different approach (see Section~\ref{sec8.1} for details).

The spectrum~$\LL_D$ was studied in~\cite{jdcs, weiss}. The
spectrum~$\LL_E$ was studied in~\cite{chaos, jdcs} (it was
introduced in~\cite{EckmannProcaccia}). See also \cite{PoW2} for
the study of transformations of the interval with an infinite
number of branches.

\section{Irregular sets and multifractal rigidity}\label{sec7}

\subsection{Multifractal analysis and irregular sets}\label{sec7.1}

In the last section we described the main components of
multifractal analysis for several multifractal spectra. These
spectra are obtained from decompositions as that in~\eqref{kkk}
and more generally as that in \eqref{kkk2}.  In particular we
possess a very detailed information from the ergodic, topological,
and dimensional point of view about the level sets in each
multifractal decomposition. However, we gave no information about
the ``irregular'' set in these decompositions, i.e., the set
$K(\varphi)$ in \eqref{kkk} and the set $X\setminus Y$ in
\eqref{kkk2}.

For example, when $\varphi\colon X\to \RR$ is a continuous
function, which is thus in $L^1(X,\mu)$ for any finite (invariant)
measure $\mu$ on~$X$, it follows from Birkhoff's ergodic theorem
that the set $K(\varphi)$ in \eqref{kkk} has zero measure with
respect to \emph{any} $S$-invariant finite measure on~$X$.
Therefore, at least from the point of view of measure theory, the
set $K(\varphi)$ is very small. However, we will see that,
remarkably, from the point of view of dimension theory, this set
is as large as the whole space, revealing once more a considerable
complexity (now for the ``irregular'' part of the multifractal
decomposition).

We first make a little digression about the concept of cohomology
in dynamical systems.  Let $S\colon X\to X$ be a continuous
transformation of the topological space~$X$.  For each function
$\varphi\colon X\to \RR$ we consider the \emph{irregular set}
$K(\varphi)$ in \eqref{kkk}. Two continuous functions
$\varphi_1\colon X\to\RR$ and $\varphi_2\colon X\to\RR$ are said
to be \emph{cohomologous} if there exists a continuous function
$\psi\colon X\to \RR$ and a constant $c\in\RR$ such that
\[
\varphi_1-\varphi_2=\psi-\psi\circ S+c
\]
on~$X$. It is easy to verify that if $\varphi_1$ and $\varphi_2$
are cohomologous then $K(\varphi_1)=K(\varphi_2)$ and
$c=P(\varphi_1)-P(\varphi_2)$. In particular, if the function
$\varphi$ is cohomologous to a constant then
$K(\varphi)=\varnothing$. The following result of Barreira and
Schmeling in \cite{ijm} shows that if $\varphi$ is not
cohomologous to a constant then $K(\varphi)$ possesses a
considerable complexity from the point of view of entropy and
Hausdorff dimension. Recall that $h(f|X)$ denotes the topological
entropy of $f|X$ (see Section~\ref{sec:for}).

\begin{theorem}[Irregular sets]\label{pijm}
If $X$ is a repeller of a $C^{1+\alpha}$ transformation, for some
$\alpha>0$, such that $f$ is conformal and topologically mixing on
$X$, and $\varphi\colon X\to\RR$ is a H\"older continuous
function, then the following properties are equivalent:
\begin{enumerate}
\item\label{apijm} $\varphi$ is not cohomologous to a constant;
\item $K(\varphi)$ is a nonempty dense set in~$X$ with
\begin{equation}\label{couve}
h(f|K(\varphi))=h(\sigma|X) \quad\text{and}\quad
\dim_HK(\varphi)=\dim_HX.
\end{equation}
\end{enumerate}
\end{theorem}

For topological Markov chains, the first identity in \eqref{couve}
was extended by Fan, Feng and Wu \cite{fan} to arbitrary
continuous functions. We note that in this case the first and
second identities in \eqref{couve} are equivalent. See also
\cite{ijm} for an appropriate version of Theorem~\ref{pijm} in the
case of hyperbolic sets.

We recall that under the hypotheses in Theorem~\ref{pijm} the set
$K(\varphi)$ has zero measure with respect to \emph{any} invariant
measure (in particular $K(\varphi)\ne X$). Theorem~\ref{pijm}
provides a necessary and sufficient condition for the set
$K(\varphi)$ to be as large as the whole space from the point of
view of entropy and Hausdorff dimension. Of~course that a priori
property~\ref{apijm} in Theorem~\ref{pijm} could be rare. However,
precisely the opposite happens.  Let $C^\theta(X)$ be the space of
H\"older continuous functions on~$X$ with H\"older exponent
$\theta\in(0,1]$ equipped with the norm
\[
\begin{split}
\lVert\varphi\rVert_\theta =&\sup\{\lvert\varphi(x)\rvert:x\in X\}\\
&+ \sup\left\{C>0:\text{$\frac{\lvert\varphi(x)-\varphi(y)\rvert}{
d(x,y)^\theta}\le C$ for each~$x$, $y\in X$}\right\}.
\end{split}
\]
It is shown in \cite{ijm} that for each $\theta\in(0,1]$ the
family of functions in $C^\theta(X)$ that are not cohomologous to
a constant forms an open dense set. Therefore, given
$\theta\in(0,1]$ and a generic function $\varphi$ in $C^\theta(X)$
the set $K(\varphi)$ is nonempty, dense in~$X$, and satisfies the
identities in \eqref{couve}.

Let now $K=\bigcup_\varphi K(\varphi)$ where the union is taken
over all H\"older continuous functions $\varphi\colon X\to\RR$.
Under the hypotheses of Theorem~\ref{pijm} we immediately conclude
that
\[
h(\sigma|K)=h(\sigma|X) \quad\text{and}\quad \dim_HK=\dim_HX.
\]
These identities were established by Pesin and Pitskel in
\cite{PP} when $\sigma$ is a Bernoulli shift with two symbols,
i.e., when
$A=\left(\begin{smallmatrix}1&1\\1&1\end{smallmatrix}\right)$ is
the transition matrix. Their methods are different from those in
\cite{ijm}. Until now it was impossible to extend the approach in
\cite{PP} even to the Bernoulli shift with three symbols.

A related result of Shereshevsky in \cite{Shereshevsky} shows that
for a generic $C^2$ surface diffeomorphism with a  locally maximal
compact hyperbolic set $\Lambda$, and an equilibrium measure $\mu$
of a H\"older continuous generic function in the $C^0$ topology,
the set
\[
I=\left\{x\in\Lambda:\liminf_{r\to0}\frac{\log\mu(B(x,r))}{\log r}
<\limsup_{r\to0}\frac{\log\mu(B(x,r))}{\log r}\right\}
\]
has positive Hausdorff dimension. This result is a particular case
of results in \cite{ijm} showing that in fact
$\dim_HI=\dim_H\Lambda$ (under those generic assumptions).

\subsection{Multifractal classification of dynamical systems}\label{sec8.2}

The former sections illustrate an enormous complexity that occurs
in a natural way in the study of the multifractal properties of
dynamical systems.  On the other hand, in the ``experimental''
study of dynamical systems it is common to have only partial
information. The multifractal spectra present themselves as
``observable'' quantities and can be determined within fairly
arbitrary precision at the expense of the macroscopic observation
of the phase space. It is thus of interest to investigate how to
recover partially or even fully the information about of a given
dynamical system, using the information contained in the
multifractal spectra. This problem belongs to the theory of
\emph{multifractal rigidity}.  In this section we want to
illustrate with a simple example how it is possible to make this
approach rigorous.

\begin{figure}[htbp]
\begin{psfrags}
\psfrag{0}{$0$} \psfrag{1}{$1$} \psfrag{A}{$A$} \psfrag{B}{$B$}
\begin{center}
\includegraphics{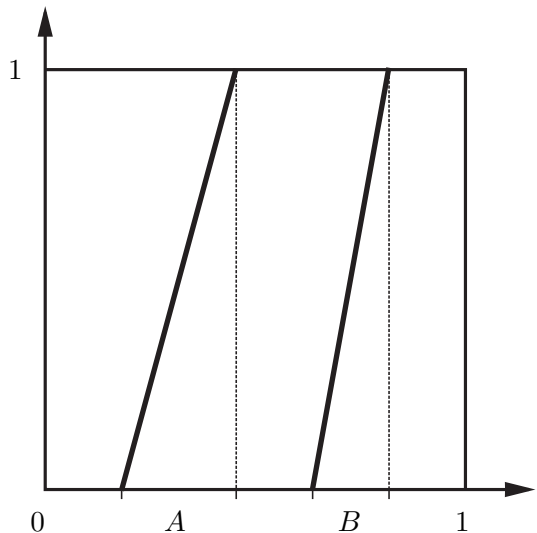}
\end{center}
\end{psfrags}
\caption{Piecewise linear expanding transformation}\label{fig8}
\end{figure}

Let $g$ and $h$ be piecewise linear transformations of the interval
$[0,1]$ with repellers given by
\[
J_g=\bigcap_{n=0}^\infty g^{-n}(A_g\cup B_g) \quad\text{and}\quad
J_h=\bigcap_{n=0}^\infty h^{-n}(A_h\cup B_h),
\]
where $A_g$, $B_g$, $A_h$ and $B_h$ are closed intervals in
$[0,1]$ such that
\[
g(A_g)=g(B_g)=h(A_h)=h(B_h)=[0,1] \quad\text{and}\quad A_g\cap
B_g=A_h\cap B_h=\varnothing.
\]
See Figure~\ref{fig8}. Both repellers can be coded by a Bernoulli
shift with two symbols.  We consider two Bernoulli measures
$\mu_g$ and $\mu_h$ (each with two symbols) invariant respectively
under $g$ and~$h$.

\begin{theorem}[Multifractal rigidity]\label{t:chaos}
If $\DD_D^{(\mu_g)}=\DD_D^{(\mu_h)}$ then there exists a
ho\-meomorphism $\chi\colon J_g\to J_h$ such that $dg=dh\circ\chi$
and $\mu_g=\mu_h\circ\chi$.
\end{theorem}

Theorem~\ref{t:chaos} is due to Barreira, Pesin and Schmeling
\cite{chaos} and provides a multifractal classification based on
the spectrum $\DD_D$.  The coincidence of the spectra of $g$ and
$h$ guarantees in particular that the derivatives of $g$ and $h$
are equal at corresponding points of the repellers and namely at
the periodic points.  There is also a version of this result for
hyperbolic sets~\cite{jdcs}. A local version of
Theorem~\ref{t:chaos} was obtained in \cite{var} for a more
general class of dynamical systems and arbitrary equilibrium
measures.

For more complex dynamical systems it may be necessary to use more
than one multifractal spectrum in order to obtain a multifractal
classification analogous to that given by Theorem~\ref{t:chaos}.
This is one of the motivations for the study of other multifractal
spectra and namely of the mixed spectra (see
Sections~\ref{sec:general} and~\ref{sec8.1}).

\section{Variational principles and number theory}\label{sec8}

\subsection{Variational principles and dimension theory}\label{sec8.1}

As we mentioned in Section~\ref{sec:general}, one can consider
several other multifractal spectra and in particular the mixed
spectra $\DD_E$ and $\EE_D$.  These two spectra combine local and
global characteristics of distinct nature, which depend not only
on the dynamics but also on the local structure provided by a
given invariant measure. As we described above (see
Section~\ref{sec6}), the spectra $\DD_D$ and $\EE_E$ are analytic
in several situations. Furthermore, they coincide with the
Legendre transform of certain functions, defined in terms of the
topological pressure. In particular, this allows one to show that
they are always convex.

In order to explain why the study of the mixed multifractal
spectra is different from the study of the nonmixed spectra, we
recall the level sets
\[
K_\alpha^{g_D}=\{x\in X:g_D(x)=\alpha\} \quad\text{and}\quad
K_\alpha^{g_E}=\{x\in X:g_E(x)=\alpha\}
\]
(see Section~\ref{sec:general}).  The main difficulty when we
study the mixed spectra is that these two families of sets need
not satisfy any Fubini type decomposition.

In the case of conformal repellers, Barreira and Saussol
\cite{mixed} obtained the following characterization of the mixed
spectra $\DD_E$ and $\EE_D$.

\begin{theorem}[Characterization of the mixed spectra]\label{lolo}
For a repeller~$X$ of a $C^{1+\varepsilon}$ transformation $f$,
for some $\varepsilon>0$, such that $f$ is conformal and
topologically mixing on~$X$, if $\varphi\colon X\to\RR$ is a
H\"older continuous function with $P(\varphi)=0$ and $\mu$ is an
equilibrium measure of $\varphi$ then
\begin{equation}\label{caca1}
\DD_E^{(\mu)}(\alpha)=\max\left\{ h_\nu(f):\text{$\nu$ is ergodic
and $-\frac{\int_X\varphi \,d\nu}{\int_X\log\lVert
df\rVert\,d\nu}=\alpha$}\right\}
\end{equation}
and
\begin{equation}\label{caca2}
\EE_D^{(\mu)}(\alpha)=\max\left\{\dim_H\nu:\text{$\nu$ is ergodic
and $-\int_X\varphi\,d\nu=\alpha$}\right\}.
\end{equation}
\end{theorem}

We call \emph{conditional variational principle} to each of the
identities in \eqref{caca1} and \eqref{caca2}. We can also obtain
conditional variational principles for the spectra $\DD_D$
and~$\EE_E$ (although in this case the results are essentially
equivalent to the corresponding multifractal analysis described in
the former sections):
\begin{equation}\label{caca3}
\DD_D^{(\mu)}(\alpha)=\max\left\{\dim_H\nu:\text{$\nu$ is ergodic
and $-\frac{\int_X\varphi \,d\nu}{\int_X\log\lVert
df\rVert\,d\nu}=\alpha$}\right\}
\end{equation}
and
\begin{equation}\label{caca4}
\EE_E^{(\mu)}(\alpha)=\max\left\{ h_\nu(f):\text{$\nu$ is ergodic
and $-\int_X\varphi\,d\nu=\alpha$}\right\}.
\end{equation}

Some of the consequences that can be obtained from Theorem~\ref{lolo}
(see \cite{mixed} for details) are the following:
\begin{enumerate}
\item the functions $\DD_E$ and $\EE_D$ are analytic; \item the
functions $\DD_E$ and $\EE_D$ are in general not convex and thus
cannot be expressed as Legendre transforms.
\end{enumerate}
The last statement strongly contrasts with what happens with the
nonmixed spectra, which are always convex.

In the case of the full shift, the identity in \eqref{caca4} was
first established by Olivier \cite{olivier, olivier2}, for the
more general class of the so-called $g$-measures. This class,
introduced by Keane in \cite{Keane}, is composed of equilibrium
measures of a class of continuous functions that need not be
H\"older continuous. It is known that any Gibbs measure is a
$g$-measure (see \cite{olivier2} for details). We note that in the
case of the shift map the identity \eqref{caca4} is equivalent to
any of the identities \eqref{caca1}, \eqref{caca2}, and
\eqref{caca3}: simply model the full shift with $m$ symbols by the
piecewise expanding map of the interval $x\mapsto mx\pmod1$ and
observe that
\[
\dim_HA=\frac{h(\sigma|A)}{\log m} \quad\text{and}\quad
\dim_H\nu=\frac{h_\nu(\sigma)}{\log m}
\]
for any ergodic $\sigma$-invariant measure $\nu$.

In \cite{mixed} the authors obtained conditional variational
principles in the more general case of transformations with upper
semi-continuous entropy (i.e., transformations $f$ for which
$\nu\mapsto h_\nu(f)$ is upper semi-continuous), for functions
with a unique equilibrium measure (and thus for a dense family of
functions; see the book by Ruelle \cite{rbb}). In \cite{tv},
Takens and Verbitski established a conditional variational
principle for the spectrum $\EE_E$ for equilibrium measures that
are not necessarily unique.

For example, when $f\colon X\to X$ is a one-sided or two-sided
topological Markov that is topologically mixing, or is an
expansive homeomorphism that satisfies specification, the entropy
is upper semi-continuous. If, in addition, $\varphi$ is a
continuous function with a certain ``bounded variation'' then it
has a unique equilibrium measure; see \cite{KH, Ke} for details.
On the other hand, one can exhibit plenty transformations that do
not satisfy specification but for which the entropy is still upper
semi-continuous. For example, all $\beta$-shifts are expansive and
thus the metric entropy is upper semi-continuous (see \cite{Ke}
for details), but in \cite{Sc} Schmeling showed that for $\beta$
in a residual set with full Lebesgue measure (although the
complement has full Hausdorff dimension) the corresponding
$\beta$-shift does not satisfy specification. It follows from work
of Walters \cite{walters2} that for every $\beta$-shift each
Lipschitz function has a unique equilibrium measure.

In \cite{waldyr}, Barreira and Saussol obtained conditional
variational principles for hyperbolic flows.

\subsection{Extensions and applications to number theory}\label{hjj}

In the theory of dynamical systems we are frequently interested in
more that one local quantity at the same time. Examples include
the Lyapunov exponents, the local entropy, and the pointwise
dimension. However, the theory described above (in
Sections~\ref{sec6} and \ref{sec8.1}) only allows us to consider
separately each of these characteristics. This observation is a
motivation to develop a multi-dimensional version of multifractal
analysis. More precisely, we want to consider, for example,
intersections of level sets of Birkhoff averages, such~as
\[
K_{\alpha,\beta}=K_\alpha(\varphi)\cap K_\beta(\psi),
\]
and to describe their multifractal properties, including their
``size'' in terms of the topological entropy and of the Hausdorff
dimension.

The corresponding multi-dimensional multifractal spectra exhibit
several nontrivial phenomena that are absent in the
one-dimensional case. Furthermore, the known approaches to the
study of one-dimensional multifractal spectra have to be modified
to treat the new situation. Nevertheless, the unifying theme
continues to be the use of the thermodynamic formalism.

We now illustrate the results with a rigorous statement in the
case of topological Markov chains. Let $\MM$ be the family of
$\sigma$-invariant probability measures on~$X$ and consider the
set
\[
\DD=\left\{\left(\int_X\varphi\,d\mu,\int_X
\psi\,d\mu\right)\in\RR^2:\mu\in\MM\right\}.
\]
The following result is a conditional variational principle for
the sets $K_{\alpha,\beta}$.

\begin{theorem}[Conditional variational principle]\label{f0}
Let $\sigma|X$ be a topologically mixing topological Markov chain,
and $\varphi$ and $\psi$ H\"older continuous functions on~$X$.
Then, for each $(\alpha,\beta)\in\interior \DD$ we have
$K_{\alpha,\beta}\ne\varnothing$ and
\begin{equation}\label{ffw}
\begin{split}
h(\sigma|K_{\alpha,\beta}) &=\sup\left\{
h_\mu(\sigma):\text{$\left(\int_X\varphi\,d\mu,
\int_X\psi\,d\mu\right)=(\alpha,\beta)$ with $\mu\in\MM$}\right\}\\
&=\inf
\left\{P(p(\varphi-\alpha)+q(\psi-\beta)):(p,q)\in\RR^2\right\}.
\end{split}
\end{equation}
\end{theorem}

Theorem~\ref{f0} is a particular case of results of Barreira,
Saussol and Schmeling~\cite{cross}. Namely, they also consider the
intersection of any finite number of level sets of Birkhoff
averages, as well as of other local quantities such as pointwise
dimensions, local entropies, and Lyapunov exponents. The first
identity in \eqref{ffw} was obtained independently by Fan, Feng
and Wu \cite{fan}, also in the more general case of arbitrary
continuous functions. In \cite{tv}, Takens and Verbitski provided
generalizations of these results. We will see below that the
second identity in \eqref{ffw} can be applied with success to
several problems in number theory.

It is also shown in \cite{cross} that if $\sigma|X$ is a
topologically mixing topological Markov chain and $\varphi$ and
$\psi$ are H\"older continuous functions on~$X$, then the
following properties hold:
\begin{enumerate}
\item if $(\alpha,\beta)\not\in\overline\DD$ then
$K_{\alpha,\beta}=\varnothing$; \item if for every $(p,q)\in\RR^2$
the function $p\varphi+q\psi$ is not cohomologous to a constant
then $\DD=\overline{\interior \DD}$; \item the function
$(\alpha,\beta)\mapsto h(\sigma|K_{\alpha,\beta})$ is analytic on
$\interior\DD$; \item there is an ergodic equilibrium measure
$\mu_{\alpha,\beta}\in\MM$ with $\int_X\varphi\,d\mu=\alpha$ and
$\int_X\psi\,d\mu=\beta$, such that
\[
\mu_{\alpha,\beta}(K_{\alpha,\beta})=1 \quad\text{and}\quad
h_{\mu_{\alpha,\beta}}(\sigma)=h(\sigma|K_{\alpha,\beta}).
\]
\end{enumerate}
In particular, the second property provides a condition which
guarantees that the identities in Theorem~\ref{f0} are valid for
an open and dense set of pairs $(\alpha,\beta)\in\DD$.

Besides their own interest, these results have several
applications and namely applications to number theory. Instead of
formulating general statements here we will describe explicit
examples that illustrate well the nature of the results obtained
by Barreira, Saussol and Schmeling in~\cite{number}. Given an
integer $m>1$, for each number $x\in[0,1]$ we denote by
$x=0.x_1x_2\cdots$ the base-$m$ representation of~$x$. It is
immediate that this representation is unique except for a
countable set of points. Since countable sets have zero Hausdorff
dimension, the nonuniqueness of the representation does not affect
the study of the dimensional properties.

For each $k\in\{0,\ldots,m-1\}$, $x\in[0,1]$ and $n\in\NN$ we define
\[
\tau_k (x,n)=\card\{i\in\{1,\ldots,n\}:x_i=k\}.
\]
Whenever there exists the limit
\[
\tau_k(x)=\lim_{n\to\infty}\frac{\tau_k(x,n)}{n}
\]
it is called the \emph{frequency} of the number $k$ in the
base-$m$ representation of~$x$.  A~classical result of Borel
\cite{Bo} says that for Lebesgue-almost every $x\in[0,1]$ we have
$\tau_k(x)=1/m$ for every~$k$. Furthermore, for $m=2$, Hardy and
Littlewood~\cite{HL} showed that for Lebesgue-almost every
$x\in[0,1]$, $k=0$,~$1$, and all sufficiently large~$n$,
\[
\left\lvert\frac{\tau_k(x,n)}n-\frac12\right\rvert
<\sqrt{\frac{\log n}{n}}.
\]
In particular, Lebesgue-almost all numbers are normal in every
integer base.  This remarkable result (that today is an immediate
consequence of Birkhoff's ergodic theorem) does not imply that the
set of numbers for which this does not happen is empty.

Consider now the set
\[
F_m(\alpha_0,\ldots,\alpha_{m-1})=\left\{x\in[0,1]
:\text{$\tau_k(x)=\alpha_k$ for $k=0$, $\ldots$, $m-1$}\right\},
\]
whenever $\alpha_0+\cdots+\alpha_{m-1}=1$ with $\alpha_i\in[0,1]$
for each~$i$. It is composed of the numbers in $[0,1]$ having a
ratio $\alpha_{k}$ of digits equal to $k$ in its base-$m$
representation for each~$k$. A precursor result concerning the
size of these sets from the point of view of dimension theory is
due to Besicovitch \cite{Be}. For $m=2$, he showed that if
$\alpha\in(0,\frac12)$ then
\[
\dim_H\left\{x\in[0,1]:\limsup_{n\to\infty}
\frac{\tau_1(x,n)}{n}\le\alpha\right\}=
-\frac{\alpha\log\alpha+(1-\alpha)\log(1-\alpha)}{\log2}.
\]
More detailed information was later obtained by Eggleston
\cite{eg}, who showed that
\begin{equation}\label{jj}
\dim_HF_m(\alpha_0,\ldots,\alpha_{m-1})=-\sum_{k=0}^{m-1}\alpha_k\log_m\alpha_k.
\end{equation}
We note that it is easy to show---and this does not require the
above result---that each set $F_m(\alpha_0,\ldots,\alpha_{m-1})$
is dense in $[0,1]$. The identity in \eqref{jj} can be established
by applying Theorem~\ref{batata} when $m=2$, Theorem~\ref{f0} when
$m=3$, and an appropriate generalization of Theorem~\ref{f0} when
$m\ge4$, thus providing a new proof of Eggleston's result (see
\cite{cross} for details).

We now consider sets of more complicated nature. Let $m=3$ and
consider the set
\[
F=\{x\in[0,1]:\tau_1(x)=5\tau_0(x)\}.
\]
This is the set of numbers in $[0,1]$ for which the base-$3$
representation has a ratio of ones that is five times the ratio of
zeros. The ratio of the digit two is arbitrary. It follows from
work in \cite{number} that
\begin{equation}\label{pp}
\dim_HF=\frac{\log(1+6/5^{5/6})}{\log3}\approx 0.85889\cdots.
\end{equation}
In order to explain how this result is obtained, we first observe
that
\begin{equation}\label{ff}
F=\bigcup_{\alpha\in[0,1/6]}F_3(\alpha,5\alpha,1-6\alpha).
\end{equation}
It is easy to show that the constant in \eqref{pp} is a lower
estimate for $\dim_HF$. Namely, it follows from \eqref{jj} and
\eqref{ff} that, since $F\supset F_3(\alpha,5\alpha,1-6\alpha)$
for each $\alpha$,
\begin{equation}\label{pp2}
\dim_HF\ge\max_{\alpha\in[0,1/6]}
-\frac{\alpha\log\alpha+5\alpha\log(5\alpha)
+(1-6\alpha)\log(1-6\alpha)}{\log3}.
\end{equation}
The maximum in \eqref{pp2} is attained at $\alpha=1/(5^{5/6}+6)$
and it is easy to verify that it is equal to the constant in
\eqref{pp}. This establishes a lower estimate for the Hausdorff
dimension.

The corresponding upper estimate is more delicate, namely because
the union in \eqref{ff} is composed of an \emph{uncountable}
number of pairwise disjoint sets.  The approach in \cite{number}
uses a generalization of the conditional variational principle in
Theorem~\ref{f0}, now for quotients of Birkhoff averages (see
\cite{number} for details).  In the particular case considered
here this variational principle states that for each $k\ne\ell$
and $\beta\ge0$ we have
\[
\dim_H \left\{x\in[0,1]:\frac{\tau_k(x)}{\tau_\ell(x)}=\beta\right\}=
\max\left\{-\sum_{j=0}^{m-1}\alpha_j\log_m\alpha_j:
\frac{\alpha_k}{\alpha_\ell}=\beta\right\}.
\]
We conclude that for each $k\ne\ell$ and $\beta\ge0$ there exists a
set
\[
F_m(\alpha_0,\ldots,\alpha_{m-1})\subset\left
\{x\in[0,1]:\frac{\tau_k(x)}{\tau_\ell(x)}=\beta\right\}
\]
with
\[
\dim_H F_m(\alpha_0,\ldots,\alpha_{m-1}) =\dim_H
\left\{x\in[0,1]:\frac{\tau_k(x)}{\tau_\ell(x)}=\beta\right\}.
\]
In particular, letting $m=3$, $k=1$, $\ell=0$, and $\beta=5$ we
conclude from \eqref{ff} that the inequality in \eqref{pp2} is in fact
an identity, and we establish \eqref{pp}.

These applications to number theory are special cases of results
established in \cite{number}: these include the study of sets
defined in terms of relations between the numbers $\tau_0(x)$,
$\ldots$, $\tau_{m-1}(x)$, and the study of sets defined in terms
of frequencies of blocks of digits, or even for which some blocks
are forbidden (thus generalizing work of Billingsley \cite{Bi}).

\section{Quantitative recurrence and dimension theory}\label{sec9}

\subsection{Quantitative recurrence}\label{secc9.1}

The Poincar\'e recurrence theorem (Theorem~\ref{thm:pp}), as
described in Section~\ref{sec1.2}, is one of the basic but
fundamental results of the theory of dynamical systems.
Unfortunately it only provides information of qualitative nature.
In particular it does not consider, for example, any of the
following natural problems:
\begin{enumerate}
\item with which frequency the orbit of a point visits a given set
of positive measure; \item\label{hju} with which rate the orbit of
a point returns to an arbitrarily small neighborhood of the
initial point.
\end{enumerate}
Birkhoff's ergodic theorem gives a complete answer to the first
problem. The second problem experienced a growing interest during
the last decade, also in connection with other fields, including
compression algorithms, and numerical study of dynamical systems.

We consider a transformation $f\colon M\to M$. The \emph{return time}
of a point $x\in M$ to the ball $B(x,r)$ (with respect to $f$) is
given by
\[
\tau_r(x)=\inf\{n\in\NN:d(f^nx,x)<r\}.
\]
The \emph{lower} and \emph{upper recurrence rates} of~$x$ (with
respect to $f$) are defined by
\begin{equation}\label{cho}
\underline R(x) =\liminf_{r\to0} \frac{\log\tau_r(x)}{-\log r}
\quad\text{and}\quad \overline R(x) =\limsup_{r\to0}
\frac{\log\tau_r(x)}{-\log r}.
\end{equation}
Whenever $\underline R(x)=\overline R(x)$ we denote the common
value by $R(x)$ and call it the \emph{recurrence rate} of~$x$
(with respect to $f$).

In the present context, the study of the quantitative behavior of
recurrence started with the work of Ornstein and Weiss \cite{ow},
closely followed by the work of Boshernitzan \cite{bosher}.
In~\cite{ow} the authors considered the case of symbolic dynamics
(and thus the corresponding symbolic metric in \eqref{g2}) and an
ergodic $\sigma$-invariant measure $\mu$, and showed that
$R(x)=h_\mu(\sigma)$ for $\mu$-almost every~$x$ (see
Section~\ref{S:produ} for more details). On the other hand,
Boshernitzan considered an arbitrary metric space $M$ and showed
in \cite{bosher} that
\begin{equation}\label{cds}
\underline R(x)\le\dim_H\mu
\end{equation}
for $\mu$-almost every~$x\in M$ (although the result in
\cite{bosher} is formulated differently, it is shown in
\cite{dimh} that it can be rephrased in this manner). It is shown
in \cite{dimh} that the inequality \eqref{cds} may be strict.

In the case of hyperbolic sets, the following result of Barreira
and Saussol in~\cite{dimh} shows that \eqref{cds} often becomes an
identity.

\begin{theorem}[Quantitative recurrence]\label{trec}
For a $C^{1+\alpha}$ diffeomorphism with a hyperbolic set
$\Lambda$, for some $\alpha>0$, if $\mu$ is an ergodic equilibrium
measure of a H\"older continuous function then
\begin{equation}\label{mnt}
R(x) =\lim_{r\to 0}\frac{\log\mu(B(x,r))}{\log r}
\end{equation}
for $\mu$-almost every point $x\in\Lambda$.
\end{theorem}

Theorem~\ref{trec} is a version of the result of Ornstein and
Weiss in \cite{ow} in the special case of symbolic dynamics. The
proof of Theorem~\ref{trec} combines new ideas with the study of
hyperbolic measures by Barreira, Pesin and Schmeling \cite{annals}
(see Section~\ref{secff}) and results and ideas of Saussol,
Troubetzkoy and Vaienti \cite{stv} and of Schmeling and
Troubetzkoy \cite{str2} (see also \cite{strextra}).
In~\cite{product}, Barreira and Saussol established a related
result in the case of repellers.

We note that the identity \eqref{mnt} relates two quantities of
very different nature.  In particular, only $R(x)$ depends on the
diffeomorphism and only the pointwise dimension depends on the
measure.

Theorem~\ref{trec} provides quantitative information about the
recurrence in hyperbolic sets. Putting together \eqref{cho} and
\eqref{mnt} we obtain
\begin{equation}\label{eq:tog}
\lim_{r\to0} \frac{\log\inf\{n\in\NN:d(f^nx,x)<r\}}{-\log r}
=
\lim_{r\to 0}\frac{\log\mu(B(x,r))}{\log r}
\end{equation}
for $\mu$-almost every point $x\in\Lambda$. Thinking as if we
could erase the limits in \eqref{eq:tog}, we can say that
Theorem~\ref{trec} shows that
\[
\inf\{k\in\NN: f^kx\in B(x,r)\}\text{ is approximately equal to
}1/\mu(B(x,r))
\]
when $r$ is sufficiently small, that is, the time that the orbit
of~$x$ takes to return to the ball $B(x,r)$ is approximately equal
to $1/\mu(B(x,r))$. This should be compared to Kac's lemma: since
$\mu$ is ergodic it tells us that
\[
\int_{B(x,r)}\tau_r(y)\,d\mu(y)=1.
\]
Hence, the average value of $\tau_r$ on $B(x,r)$ is equal to
$1/\mu(B(x,r))$. Therefore, Theorem~\ref{trec} can be though of as
a local version of Kac's lemma.

In another direction, the results in \cite{dimh, product} motivate
the introduction of a new method to compute the Hausdorff
dimension of a given measure (see Section~\ref{sec:dd}). See
\cite{dimh} for details.

We now consider the case of repellers and briefly present two
applications to number theory. Let $x=0.x_1x_2\cdots$ be the
base-$m$ representation of the point $x\in[0,1]$.  It was shown
in~\cite{product} that
\[
\text{$\inf\{n\in\NN\colon \lvert
0.x_nx_{n+1}\cdots-0.x_1x_2\cdots\rvert<r\}\sim\frac1r$ when
$r\to0$}
\]
for Lebesgue-almost every $x\in[0,1]$, meaning that
\[
\lim_{r\to0} \frac{\log\inf\{n\in\NN\colon \lvert
0.x_nx_{n+1}\cdots-0.x_1x_2\cdots\rvert<r\}}{-\log r}=1
\]
for Lebesgue-almost every $x\in[0,1]$.  Another example is given
by the continued fractions. Writing each number $x\in(0,1)$ as a
continued fraction
\[
x=[m_1,m_2,m_3,\ldots]
=\cfrac{1}{m_1+\cfrac{1}{m_2+\cfrac{1}{m_3+\cdots}}},
\]
with $m_i=m_i(x)\in\NN$ for each $i$ (this representation is unique
except for a countable subset of $(0,1)$), it is shown in
\cite{product} that
\[
\text{$\inf\{n\in\NN:\lvert [m_n,m_{n+1},\ldots]-[m_1,m_2,\ldots]
\rvert<r\}\sim\frac1r$ when $r\to0$}
\]
for Lebesgue-almost every $x\in(0,1)$.

\subsection{Product structure and recurrence}\label{S:produ}

We already described the product structure of hyperbolic sets (see
Section~\ref{sec1.1}) and the product structure of hyperbolic
measures (see Section~\ref{secff}). The study of quantitative
recurrence can also be used to obtain new information about the
product structure.

\begin{figure}[htbp]
\begin{psfrags}
\psfrag{FFF}{$[f^nx,x]$} \psfrag{Vxx}{$V^u_\varepsilon(x)$}
\psfrag{VV}{$V^s_\varepsilon(x)$} \psfrag{X}{$x$}
\psfrag{Fnx}{$f^nx$} \psfrag{BBB}{$B(x,\rho)$}
\begin{center}
\includegraphics{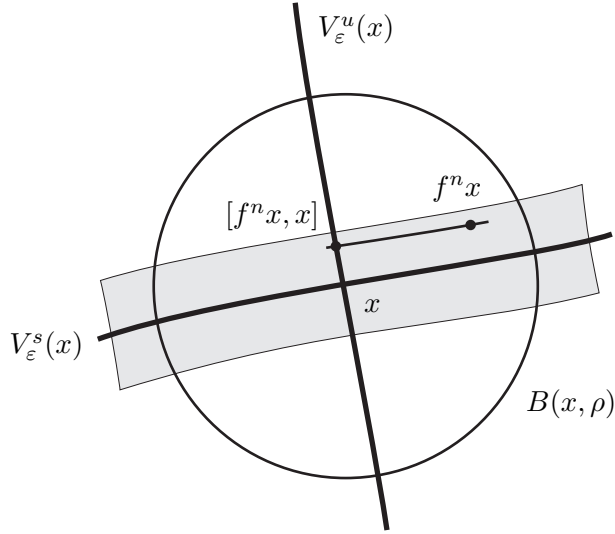}
\end{center}
\end{psfrags}
\caption{Definition of the unstable return time (the shaded area
is the set of points at a $d^u$-distance of $V^s_\varepsilon(x)$
at most $r$)}\label{fig9}
\end{figure}

Let $f\colon M\to M$ be a $C^{1+\alpha}$ diffeomorphism with a
locally maximal compact hyperbolic set $\Lambda\subset M$, and
denote by $d_s$ and $d_u$ the distances induced by the distance
$d$ of $M$ respectively on each stable and unstable manifold. When
$d(f^nx,x)\le\delta$, for each $\rho\le\delta$ we can define (see
Figure~\ref{fig9})
\[
\begin{split}
\tau_r^s(x,\rho) &= \inf\{n\in\NN:\text{$d(f^{-n}x,x) \le \rho$
and $d_s([x,f^{-n}x],x)<r$}\},\\ \tau_r^u(x,\rho) &=
\inf\{n\in\NN:\text{$d(f^nx,x) \le \rho$ and
$d_u([f^nx,x],x)<r$}\}.
\end{split}
\]
We call $\tau_r^s(x,\rho)$ and $\tau_r^u(x,\rho)$ respectively
\emph{stable} and \emph{unstable return times}.  Note that the
functions $\rho\mapsto\tau_r^s(x,\rho)$ and
$\rho\mapsto\tau_r^u(x,\rho)$ are nondecreasing.  We~define the
\emph{lower} and \emph{upper stable recurrence rates} of the point
$x\in\Lambda$ (with respect to $f$) by
\[
\underline R^s(x) = \lim_{\rho\to0} \underline R^s(x,\rho)
\quad\text{and}\quad \overline R^s(x) = \lim_{\rho\to0} \overline
R^s(x,\rho),
\]
and the \emph{lower} and \emph{upper stable recurrence rates} of
the point $x\in\Lambda$ (with respect to $f$) by
\[
\underline R^u(x) = \lim_{\rho\to0} \underline R^u(x,\rho)
\quad\text{and}\quad \overline R^u(x) = \lim_{\rho\to0} \overline
R^u(x,\rho),
\]
where
\[
\underline R^s(x,\rho) = \liminf_{r\to0}
\frac{\log\tau_r^s(x,\rho)}{-\log r} \quad\text{and}\quad
\overline R^s(x,\rho) = \limsup_{r\to0}
\frac{\log\tau_r^s(x,\rho)}{-\log r},
\]
\[
\underline R^u(x,\rho) = \liminf_{r\to0}
\frac{\log\tau_r^u(x,\rho)}{-\log r} \quad\text{and}\quad
\overline R^u(x,\rho) = \limsup_{r\to0}
\frac{\log\tau_r^u(x,\rho)}{-\log r}.
\]
When $\underline R^s(x)=\overline R^s(x)$ we denote the common
value by $R^s(x)$ and call it \emph{stable recurrence rate} of~$x$
(with respect to $f$), and when $\underline R^u(x)=\overline
R^u(x)$ we denote the common value by $R^u(x)$ and call it
\emph{unstable recurrence rate} of~$x$ (with respect to~$f$).

Barreira and Saussol showed in \cite{product} that for a
$C^{1+\alpha}$ diffeomorphism that is topologically mixing on a
locally maximal compact hyperbolic set $\Lambda$, and an
equilibrium measure $\mu$ of a H\"older continuous function, we
have
\begin{equation}\label{utp}
R^s(x)=\lim_{r\to0}\frac{\log\mu^s_x(B^s(x,r))}{\log r}
\quad\text{and}\quad
R^u(x)=\lim_{r\to0}\frac{\log\mu^s_x(B^u(x,r))}{\log r}
\end{equation}
for $\mu$-almost every $x\in\Lambda$, where $\mu^s_x$ and
$\mu^u_x$ are the conditional measures induced by the measurable
partitions $\xi^s$ and $\xi^u$ (see Section~\ref{secff}).
Ledrappier and Young~\cite{LY} showed that there exist the limits
in the right-hand sides of the identities in \eqref{utp} (see
Section~\ref{sec:dd}).

The following result in \cite{product} can now be obtained using
Theorems~\ref{thm:annals} and~\ref{trec} and the identities in
\eqref{utp}.

\begin{theorem}[Product structure for recurrence]\label{thmaxa}
Let $\Lambda$ be a locally maximal compact hyperbolic set of a
$C^{1+\alpha}$ diffeomorphism that is topologically mixing
on~$\Lambda$, for some $\alpha>0$, and $\mu$ an equilibrium
measure of a H\"older continuous function. Then, for $\mu$-almost
every point $x\in \Lambda$ the following properties hold:
\begin{enumerate}
\item
the recurrence rate is equal to the sum of the stable and unstable
recurrence rates, i.e.,
\[
R(x)=R^s(x)+R^u(x);
\]
\item there exists $\rho(x)>0$ such that for each $\rho<\rho(x)$
and each $\varepsilon>0$ there is $r(x,\rho,\varepsilon)>0$ such
that if $r<r(x,\rho,\varepsilon)$ then
\[
r^\varepsilon < \frac{ \tau_r^s(x,\rho) \cdot \tau_r^u(x,\rho) }{
\tau_r(x) } < r^{-\varepsilon}.
\]
\end{enumerate}
\end{theorem}

The second statement in Theorem~\ref{thmaxa} shows that the return
time to a given set is approximately equal to the product of the
return times in the stable and unstable directions, as if they
were independent.

A related result was obtained by Ornstein and Weiss in the case of
symbolic dynamics. Namely, they showed in \cite{ow} that if
$\sigma^+\colon\Sigma^+\to\Sigma^+$ is a \emph{one-sided} subshift
and $\mu^+$ is an ergodic $\sigma^+$-invariant probability measure
on $\Sigma^+$, then
\begin{equation}\label{ow1}
\lim_{k\to\infty} \frac{\log\inf\{n\in\NN:(i_{n+1}\cdots
i_{n+k})=(i_1\cdots i_k)\}}{k}=h_{\mu^+}(\sigma)
\end{equation}
for $\mu^+$-almost every $(i_1i_2\cdots)\in\Sigma^+$. They also
showed in \cite{ow} that if $\sigma\colon\Sigma\to\Sigma$ is a
\emph{two-sided} subshift and $\mu$ is an ergodic
$\sigma$-invariant probability measure on~$\Sigma$, then
\begin{equation}\label{ow2}
\lim_{k\to\infty} \frac{\log\inf\{n\in\NN:(i_{n-k}\cdots
i_{n+k})=(i_{-k}\cdots i_k)\}}{2k+1}=h_\mu(\sigma)
\end{equation}
for $\mu$-almost every $(\cdots i_{-1}i_0i_1\cdots)\in\Sigma$.

Given a two-sided shift $\sigma\colon\Sigma\to\Sigma$ it has
naturally associated two one-sided shifts
$\sigma^+\colon\Sigma^+\to\Sigma^+$ and
$\sigma^-\colon\Sigma^-\to\Sigma^-$ (respectively related with the
future and with the past). Furthermore, any $\sigma$-invariant
measure $\mu$ on $\Sigma$ induces a $\sigma^+$-invariant measure
$\mu^+$ on $\Sigma^+$ and a $\sigma^-$-invariant measure $\mu^-$
on $\Sigma^-$, such that
\[
h_{\mu^+}(\sigma^+)=h_{\mu^-}(\sigma^-)=h_{\mu}(\sigma).
\]
For each $\omega=(\cdots i_{-1}i_0i_1\cdots)\in\Sigma$ and
$k\in\NN$ we set
\[
\begin{split}
\tau^+_k(\omega)
&=\inf\{n\in\NN:(i_{n+1}\cdots i_{n+k})=(i_1\cdots i_k)\},\\
\tau^-_k(\omega)
&=\inf\{n\in\NN:(i_{-n-k}\cdots i_{-n-1})=(i_{-k}\cdots i_{-1})\},\\
\tau_k(\omega) &=\inf\{n\in\NN:(i_{n-k}\cdots
i_{n+k})=(i_{-k}\cdots i_k)\}.
\end{split}
\]
Let now $\mu$ be an ergodic $\sigma$-invariant measure on
$\Sigma$. It follows from \eqref{ow1} and \eqref{ow2} that for
$\mu$-almost every $\omega\in\Sigma$, given $\varepsilon>0$, if
$k\in\NN$ is sufficiently large then
\[
e^{-k\varepsilon}\le\frac{\tau^+_k(\omega)
\tau^-_k(\omega)}{\tau_k(\omega)}\le e^{k\varepsilon}.
\]
Theorem~\ref{thmaxa} and the identities in \eqref{utp} are
versions of these statements in the case of dimension.

\bibliographystyle{amsplain}

\enlargethispage{.2cm}

\end{document}